\documentclass[final,3p]{elsarticle}

\usepackage{amsmath}
\usepackage{mathtools}
\usepackage{amssymb}
\usepackage{amsbsy}
\usepackage{mathrsfs}
\usepackage{graphicx}
\usepackage{pgfplots}
\usepackage{tikz}
\usepackage{multirow}

\def\Vu{ {\bf u} }
\def\pVu{ {\tilde{\bf u} }}
\def\Vv{ {\bf v} }
\def\Vn{ {\bf n} }

\def\vp{\text{p}}
\def\VF{ {\bf F} }

\def\V{ {\bf V} }
\def\VD{ {\bf D} }
\def\VS{ {\bf S} }

\DeclareMathOperator*{\argmin}{argmin}

\def\Vs#1{{\boldsymbol{#1}}}
\def\Du{\boldsymbol D(\Vs{u})}

\usepackage{xcolor}

\newtheorem{thm}{Theorem}
\newtheorem{remark}{Remark}
\newtheorem{hypothesis}{Hypothesis}

\journal{}%Applied Mathematics Letters}

\begin{document}

\begin{frontmatter}
\title{Shear rate projection schemes for non--Newtonian fluids.\tnoteref{CRSNG}}

\author[Giref]{J. ̃Deteix\corref{cor1}}
\ead{jean.deteix@mat.ulaval.ca}
\author[Giref]{D. ̃Yakoubi}
\ead{yakoubi@giref.ulaval.ca}

\address[Giref]{Groupe Interdisciplinaire de Recherche en \'El\'ements Finis de l'Universit\'e Laval, D\'epartment de Math\'ematiques et Statistiques, Universit\'e Laval,Qu\'ebec,Canada}

\cortext[cor1]{Corresponding author}
\tnotetext[CRSNG]{This work was supported by the Natural Sciences and Engineering Research Council of Canada (NSERC) Discovery Grant \# RGPIN-2015-04932 (J. Deteix).}

%\maketitle
\begin{abstract}
The operator splitting approach applied to the Navier-Stokes equations, gave rise to various numerical methods for the simulations of the dynamics of fluids. The separate work of Chorin and Temam on this subject gave birth to the so-called \textit{projection methods}. 
The basic projection schemes, either the \textit{incremental} or \textit{non-incremental} variant (see \cite{GueMinShe2006}) induces an artificial Neumann boundary condition on the pressure. By getting rid of this boundary condition on the pressure, the so-call \textit{rotational incremental pressure-correction scheme} as proposed by  Timmermans et al.~\cite{TimMinVan1996} for Newtonian fluids with constant viscosity gives a consistent equation for the pressure. In this work we propose a family of projection methods for \textit{generalized Newtonian fluids} based on an extension of the rotational projection scheme. Called \textit{shear rate projections}, these methods produces consistent pressure when applied to generalized Newtonian fluids. Accuracy of the methods will be illustrated using a manufactured solution. Numerical experiments 
%for the lid-driven cavity flow and 
for the flow past a cylinder, %both 
with a Carreau rheological model, will also be presented. 
%In the natural convection case, this new approach will be combined with the coupled projection scheme (see \cite{DetJenYak2014}) leading to a new coupled projection method applicable to generic coupled convection-diffusion and Navier--Stokes equations.
\end{abstract}

\begin{keyword}
% Ancien texte
% Navier Stokes, Equation, heat Equation, fractional time stepping,
% polynomial

% Nouveau texte
Navier--Stokes equations, non--Newtonian, Carreau, projection scheme, shear rate.

\MSC 76D05 \sep  65M60 \sep 35Q35

\end{keyword}

\end{frontmatter}
 \section{Introduction}
Solving the Navier--Stokes system  describing unsteady flows is a theoretical and numerical challenge for a homogeneous fluid and, \textit{a fortiori}, for heterogeneous fluids. To numerically simulate such flows, apart from the heterogeneous nature of the fluid,  we are confronted with two major difficulties : very complex dynamics, needing accurate time approximation (high order time approximation or fine time step); complex spatial behaviour that dictates the use of very fine meshes when solving the incompressible Navier--Stokes equations. The construction of efficient solver to achieve good approximation in a reasonable computational time is a difficult task. Since we are interested in three dimensional problems, numerical methods based on mixed formulation (for example \cite{BofBreFor2013, GirRav1986} in the finite element case), although precise, can be time and resources consuming. %\JDm{references Newtonien-Non-Newtonien SVP}

%One way to circumvent the intrinsic saddle point structure due to the
%incompressibility constraint of the Navier--Stokes equation, is to follow
In designing an effective time-marching techniques for this problem we need to address the
fact that the incompressibility constraint in the Navier--Stokes system gives the problem a saddle point structure.
One way to overcome this structure is found in the pioneering works of Chorin~\cite{ChoAle1968,ChoAle1969}  and  Temam~\cite{Tem1977}
who introduced projection methods. % which are special cases of fractional time-stepping schemes.
The idea is to apply a fractional time step to decouple the  incompressibility constraint from the diffusion operator based on the Helmoltz decomposition (see \cite{GirRav1986} for instance).

The projection method as originally presented has some drawbacks: limited precision of the resulting algorithm (see Rannacher \cite{Ran1992} and Shen \cite{She1992}, the use of any higher-order 
time stepping scheme does not improve the overall accuracy), and artificial boundary condition of the resulting pressure.

To alleviate those deficiencies, numerous variants have been proposed over the years (\cite{BelColGla1989,Ran1992,She1996,Gue1997,GueQua1998}).
For an interesting overview, we also refer to Guermond and coauthors~\cite{GueMinShe2006,GueQua1998}.
Notably, those limitations lead to the introduction of an \textit{incremental projection scheme} (proposed by Goda et al \cite{God1979,Van1986}), a projection where the viscous equation takes into account the pressure at the previous time-step. 
In the context of homogeneous fluids, this variation of the original scheme is commonly used since it provides an improved precision (see \cite{GueMinShe2006}).

%===================

Timmermans et al. \cite{TimMinVan1996} proposed the \textit{rotational projection scheme} (RP) or \textit{rotational incremental pressure-correction scheme} in Guermond's terminology.
This scheme gives a consistent boundary condition for the pressure and better precision for the velocity and pressure, see for instance Guermond et al.  \cite{GueMinShe2006}.
However, this approach is only valid for \textit{homogeneous viscosity} as illustrated in~\cite{DetYak2017}; limiting its use to homogeneous Newtonian flows.
Still based on the original idea of Chorin and Temam, the more general \textit{shear rate projection} (SRP) proposed by Deteix et al.~\cite{DetYak2017} allows the treatment of heterogeneous viscosity and improves the accuracy of the incremental projection (the rate of convergence are comparable to those of the  rotational projection for Newtonian fluids).
This enhanced precision makes the shear rate projection attractive in applications such as natural convection, Allen-Cahn or Cahn-Hilliard flows but also %in generalized Newtonian 
non--Newtonian fluids.

Non--Newtonian model can be summarily described as modelling complex flows using  
%Complex non--Newtonian fluids obey 
constitutive models which involve dependency of the local stress on the velocity gradient (shear-rate dependant fluid) and possibly the deformation history of the fluid. In this work we are interested in \textit{generalized Newtonian fluids}.  
%models where the velocity fields is instantaneously modified by the stresses, that is a the constitutive relation 
%models where stresses depends on velocity and its gradient but not on the history of the flow. 
Generalized Newtonian fluids form a subclass of non--Newtonian models. They are the simplest extension from Newtonian models to non--Newtonian models. %Sometimes called shear-rate dependant fluids, 
They model phenomenon where the flow affect the viscosity of the fluid, but does not change the "Newtonian nature" of the constitutive law :
$$\sigma = -p I + \nu(\|\Du\|)\Du.$$ 
%As an extension of the Newtonian model 
The velocity fields of a generalized Newtonian fluids is instantaneously modified by the stresses (the history of the flow has no effect). 
These fluids are characterized by the  derivative of the viscosity with respect to $\|\VD(u)\|$ (the tensorial norm).  Monotone viscosity leads to two possible states: \textit{shear-thickening} where resistance to shear increases as the shear rate increases and s\textit{hear-thinning} fluids having the opposite behaviour. Of course more complex behaviour can be modelled with non-monotone relation between shear-rate and viscosity.

The constitutive law is purely phenomenological as it try to replicate the stresses due to the applied flow by a shear-dependent effective viscosity. This give rise to various phenomenological description of viscosity adapted to various empirical observations. \textit{Carreau--Yasuda},  \textit{power law}, \textit{Cross}, etc.  \cite{OwePhi2002, Irg2013} are constitutive equations for $\nu$ of particular interest as they seems to be frequently used in numerous situations.
There is an extensive and recent literature (\cite{BaeWol2016, TanZho2017, DreHun2004, DieRuzWol2010, Bae2015} for example) concerning existence and regularity of solutions of the generalized Newtonian flow problem for multiple constitutive laws. 

%%% Deplacer dans le papier SRP+THermo-fluid
% In~\cite{DetJenYak2014}, using a method (called \textit{coupled projection}) based on the incremental projection, an increase in computational performance is observed compared to usual approaches for natural convection problems. Similar gain in performance are more than likely for coupled problems involving non Newtonain flows (especially in convection-diffusion, coextrusion for example). 
%%%

Our goal is to propose an original projection method for the numerical simulation of generalized Newtonian flow. To achieve this, we will restrict ourselves to a rheological and theoretical context insuring validity of the continuous model (more precisely we refer to \cite{DreHun2004}).  
% Unfortunately the 
The shear rate projection as proposed in \cite{DetYak2017} does not take into account the explicit dependence of the viscosity upon velocity. However, the same paradigm applies to the projection method when used on generalized Newtonian fluid:  some information, easily extracted from the predicted velocity, is ignored. We propose a modification of the shear rate projection, recuperating part of the lost information, making it possible to gain in accuracy and consistency. Based on this generalized version of the SRP, a family of numerical scheme based on the finite element method combined to a second order time approximation is proposed.
%Pas besoin d'\'evoquer les autres m\'ethodes, on propose la notre.\\

The last section will be devoted to explore different representative of this family of methods. The implementation is based on the finite element library FreeFem, \cite{freefem}, code efficiency and optimal strategies will be minimal and performance measure will be limited to comparison of physical quantities for various methods. 
% However, 
Time and spatial accuracy will be illustrated first, followed by a classical application, that can be regarded as benchmark: the steady flow past a cylinder. 
\section{Problem setting}\label{Position} %\hspace{20cm}
We consider an unsteady flow of a generalized Newtonian fluid : a fluid having a non-homogeneous viscosity $\nu$ depending on time, space, velocity or shear rate of the fluid, possibly the fluid pressure, temperature or other external quantities and a constitutive equations of the form
$$ \sigma(\Vs{u}) = -pI + 2\nu\Du.$$
Assuming a Boussinesq-type model, the density variations are neglected with the exception of the forces term. These assumptions leads to   
%we are aiming at solving 
the Navier--Stokes equations for the fluid velocity $\Vu$ and pressure $p$  
%and the temperature $T$
%, which lead to the following system
\begin{equation}\label{Prob_Boussinesq1}
\left \{
\begin{aligned}
\displaystyle \rho&\frac{\partial \Vs{u}}{\partial t}   + \rho\left(\Vs{u} \cdot \nabla \right)\Vs{u}
- \nabla \cdot(2\nu \Du) + \nabla \,p =\Vs{f},\\
&\nabla \cdot \Vs{u} = 0
\end{aligned}
\right.
\qquad\mbox{ on } \Omega_t,
% \left \{
% \begin{array}{rcll}
% \displaystyle \rho\frac{\partial \Vu}{\partial t}   + \rho\left(\Vu \cdot \nabla \right)\Vu
% - \nabla \cdot(2\nu \VD(\Vu)) + \nabla \,p &=&
% \Vs{f},  \qquad \,  \, & \\
% \nabla \cdot \Vu &=& 0,&\\
% %\displaystyle \frac{\partial T}{\partial t}   +  \left(\Vu \cdot \nabla \right) T
% %- \nabla \cdot(\lambda(T) \nabla T) &=& H& \mbox{ in } \Omega_t,
% \end{array}
% \right. \mbox{ in } \Omega_t,
\end{equation}
%where the function $H$ represents an external heat source and depends only on the position vector
where  the shear rate is defined as
$$
\Du = \frac{1}{2}(\nabla\Vs{u} + \nabla\Vs{u}^t).
$$
The general expression
$$\nu(t,\Vs{u}) = \nu(t,\Vs{x},\Vs{u}, \Du, ...)$$
represents the viscosity of the fluid and $\Vs{f}$  represents external volumic forces (such as gravity).
System \eqref{Prob_Boussinesq1} is completed with the following initial data:
\begin{equation}\label{BC-Initial}
\Vs{u}(0,\Vs{x}) = \Vs{u}_0(\Vs{x}) \in L^2(\Omega)^d \; \mbox{with} \; \nabla \cdot \Vs{u}_0 = 0 \;
\end{equation}
%     which are assumed to belong to $L^2(\Omega)^d$ and $L^2(\Omega)$ respectively
and we  consider homogeneous Dirichlet boundary conditions on a non-empty, but possibly limited, part of the boundary of the domain%for $\Vu$ and both $\Gamma_D$ and $\Gamma_N$ of positive measure
%     The homogenous Dirichlet boundary condition $\Vu = 0$ on $\partial \Omega$  can be extended
%but the general case follows the same lines. %(for example Girault-Raviart \cite{GirRav1986})
%     following the Hopf Lemma \cite{Hop1955}.
% \item
%and boundary condition
\begin{equation}\label{BC-Limit}
\Vs{u} = 0 \; \mbox{ on} \, \Gamma_D \subseteq \partial \Omega.
\end{equation}
%     \begin{remark}
% \begin{itemize}
%  \item

% Both $\Gamma_D$ and $\Gamma_N$ are of positive measure. Of course this assumption 
% could be omitted by imposing a suitable compatibility condition.
% \end{itemize}
%     \end{remark}
\subsection{Existence and regularity of solutions}
Let $\Omega$  be a smooth domain in $\mathbb R^d,~ d=$  2 or 3 which satisfy the inf--sup conditions (see \cite{GirRav1986}). Let $\partial \Omega= \Gamma_D \cup \Gamma_N$ the boundary of $\Omega$ 
(possibly $\Gamma_N = \emptyset$)), $\Vs{n}$ is the exterior normal vector and $\Omega_t$ the open set  $\Omega \times \left(0,T\right)$, 
where $T > 0$ is the final time.

We denote by $L^{2}(\Omega)$ the space of square integrable functions defined on $\Omega$.  The Sobolev spaces $W^{p,q}(\Omega)$, $q \geq 0$, are the spaces of functions in $L^{q}(\Omega)$ with generalized partial derivatives belonging to $L^{q}(\Omega)$ up to order $p$ (see \cite{Ada2003}). For $X$ a Banach space, the Bochner spaces $L^r(0,T, X)$ are the spaces of functions $v: t \mapsto v(t)$ defined on $(0,T)$ with values in $X$ (see \cite{Ada2003}).
%The Sobolev spaces are equipped with the usual norm $\|\cdot\|_{m, \Omega}$ (see \cite{adams}). 
The functional space $\mathscr{W}_r$ and $\mathscr{M}$ %containing the velocity $\Vu$ and pressure $p$ respectively 
are defined as
\begin{equation*}
\mathscr{W}_r = \left\{\Vs{v} \in L^{r}(0, T; V_r): \partial_t\Vs{v} \in L^{r/(r-1)}(0, T; V_r')\right\},\qquad
\mathscr{M} = L^{2}(0, T; L^2(\Omega))
\end{equation*}
with
$$V_r = \left\{\Vs{v} \in \left(W^{1,r}(\Omega)\right)^d: \nabla\cdot\Vs{v} = 0 \text{ on }\Omega,\,\,\, \Vs{v} = 0 \mbox{ on }  \Gamma_{D}\right\}.$$
%\newpage
%

Regarding the mathematical analysis of this model, specifically the existence and regularity of solutions, we refer the readers to works such as \cite{BaeWol2016,TanZho2017,DreHun2004} and the references therein. 
To insure existence of a solution of \eqref{Prob_Boussinesq1}-\eqref{BC-Limit} in case of non--Newtonian fluids some assumptions on $\sigma(\Vs{u})$ are needed. 
%(see \cite{DieRuzWol2010,DreHun2004} and the references therein). 
The results proposed in~\cite{DreHun2004}, based on three general assumptions, is presented here. %
Since we are interested in generalized Newtonian fluid, two of these assumptions, here noted \eqref{CondNu0} and \eqref{CondNu1}, are reformulated putting the emphasis on $\nu(t,\Vs{x})$, the viscosity of the fluid. 
Notice that rheological models such as the power law and Carreau--Yasuda type models used in a wide variety of industrial applications all verify these assumptions with little restrictions (see Remark~\ref{rmkCarYas}).

\begin{hypothesis}[continuity]\label{CondNu0} 
$$\begin{array}{c}
\nu: (0,T)\times\Omega\times \mathbb{R}^d\times \mathbb{R}^{d\times d} \longrightarrow \mathbb{R}\vspace*{4pt}\\
(t,\Vs{x},\Vs{u},\Vs{F}) \longmapsto \nu(t,\Vs{x},\Vs{u},\Vs{F})\\
\end{array}
$$
is measurable with respect to $(t,\Vs{x})$ for all $(\Vs{u},\Vs{F}) \in \mathbb{R}^d\times \mathbb{R}^{d\times d}$ and continuous with respect to $(\Vs{u},\Vs{F})$ for almost every $(t,\Vs{x})\in (0,T)\times\Omega$.
\end{hypothesis}
\begin{hypothesis}[coercivity]\label{CondNu1}
There exist $c_1\ge 0, c_2 > 0, \lambda_1 \in L^{r'}((0,T)\times\Omega), 
\lambda_2 \in L^{1}((0,T)\times\Omega), \lambda_3 \in L^{(r/\alpha)'}((0,T)\times\Omega), 0<\alpha<r$, such that
$$ \nu(t,\Vs{x},\Vs{u},\Vs{F}) \le \|\Vs{F}\|^{-1}\left(\lambda_1(t,x)+c_1\left(\|\Vs{u}\|^{r-1}+\|\Vs{F}\|^{r-1}\right)\right)$$
$$ \nu(t,\Vs{x},\Vs{u},\Vs{F}) \ge \|\Vs{F}\|^{-2}\left(-\lambda_2(t,x)-\lambda_3(t,x)\|\Vs{u}\|^{\alpha}+c_2\|\Vs{F}\|^{r}\right)$$
with $\|\VF\|^2 = \sum_{ij} \VF_{ij}^2$ the Frobenius norm of the tensor $\VF$.
\end{hypothesis}
\begin{hypothesis}[monotonicity]\label{CondNu2}
%$\nu$ satisfies one of the following conditions:

%$\bullet\,
$\forall\ (t,\Vs{x}, \Vs{u})\in (0,T)\times\Omega\times\mathbb{R}$ , the map $\Vs{F}\mapsto\nu(t,\Vs{x},\Vs{u},\Vs{F})$ is a $C^1$ function and is monotone,
$$\left(\nu(t,\Vs{x},\Vs{u},\Vs{F})\Vs{F} - \nu(t,\Vs{x},\Vs{u},\Vs{G})\Vs{G}\right):(\Vs{F}-\Vs{G}) \ge 0\qquad\forall\ \Vs{F}, \Vs{G}\in \mathbb{R}^{d\times d}.$$
% $\bullet\,\forall\ (t,\Vs{x}, \Vs{u})\in (0,T)\times\Omega\times\mathbb{R}$ , the map $\Vs{F}\mapsto\nu(t,\Vs{x},\Vs{u},\Vs{F})$ is strictly monotone,
% $$\left(\nu(t,\Vs{x},\Vs{u},\Vs{F})\Vs{F} - \nu(t,\Vs{x},\Vs{u},\Vb{G})\Vb{G}\right):(\Vs{F}-\Vb{G}) > 0\qquad \forall\ \Vs{F}\neq\Vb{G}\in \mathbb{R}^{d\times d}.$$
\end{hypothesis}
%
% \begin{hypothesis}\label{CondNu3}
% We assume that $\nu$ admit a strictly positive lower bound
% \begin{equation}\label{borne_nu}
% 0 < \nu_0   \leq \nu 
% %\qquad \forall t \in \mathbb [0,T],\, \Vs{x}\in \Omega
% \end{equation}
% \end{hypothesis}
%
\begin{thm}
Assuming $\nu$ satisfies hypothesis~\eqref{CondNu0}-\eqref{CondNu2} for some $r \in [\frac{3d+2}{d+2}, \infty)$. Then for  $\Vs{u}_0 \in \left(L^2(\Omega)\right)^d$ and $\Vs{f} \in L^{r/(r-1)}(0,T;V_r')$, the system \eqref{Prob_Boussinesq1}-\eqref{BC-Limit} has a weak solution $(\Vs{u},p)$, with $u\in \mathscr{W}_r$ and $p\in \mathscr{M}$.
\end{thm}

\begin{remark}\label{rmkCarYas}
Following \cite{DreHun2004,Bae2015,DieRuzWol2010} the conditions~\eqref{CondNu0}-\eqref{CondNu2} are verified in the case of a generic law 
\begin{equation}\label{loiCarYas}
\nu(\Du) = \nu_\infty + (\nu_0-\nu_\infty)
\displaystyle\left(C_0 + \lambda^2 \|\Du\|^2\right)^{\frac{m-1}{2}}\qquad m, C_0, \nu_\infty,\nu_0 \ge 0.%,  \nu_\infty \neq \nu_0.
\end{equation}
Thanks to the numerous parameters, this correspond to various rheological models (power law, Carreau, Cross, Carreau--Yasuda, etc.) which are frequently used in engineering context. In such cases, from \cite{DieRuzWol2010} we get existence of a weak solution, provided $m> 0$ for $d=2$ and $m> 1/5$ for $d=3$.
\end{remark}

\subsection{Time discretization}

The choice of time discretization is motivated by the fact that we rely on splitting technique to construct approximation of \eqref{Prob_Boussinesq1}. Projection schemes (as methods based on operator splitting) have an inherent splitting error of order $3/2$ in $H^1-$norm (see \cite{Ran2000,GueShe2004,GueMinShe2006}). 
Therefore the proposed algorithm, relying on a projection scheme, is \textit{at best} of second order in $H^1-$norm, and the use of higher order time discretization is irrelevant.

For the sake of simplicity and clarity, the method proposed in this paper will be based on the second order backward time discretization (BDF2) which is frequently used to solve Navier--Stokes equation \eqref{Prob_Boussinesq1}. 
$$
D_t\Vs{u}^{n+1} = \frac{{3\Vs{u}}(t^{n+1},\Vs{x})-4\Vs{u}(t^n,\Vs{x})+\Vs{u}(t^{n-1},\Vs{x})}{2\Delta t}.
$$
Using a constant time-step $\Delta t$, denoting $t^n = n\Delta t$, $\Vs{u}^{n} = \Vs{u}(t^{n}, \Vs{x})$,  $p^n = p(t^n,\Vs{x})$, $\Vs{f}^n = \Vs{f}(t^n,\Vs{x})$ and $\nu^{n}(\Vs{u}) = \nu(t^{n}, \Vs{u})$ the implicit time discretization of  \eqref{Prob_Boussinesq1} 
result in a sequence of (generalized) Oseen problems of the form
\begin{equation}\label{Prob_TimeDiscretization}
\left \{
\begin{aligned}
\rho &D_t\Vs{u}^{n+1}   + \rho \left(\Vs{u}^{n+1} \cdot \nabla \right)\Vs{u}^{n+1}  
%\\ & \qquad\qquad\qquad\qquad
-\nabla \cdot(2\nu^{n+1}(\Vs{u}^{n+1}) \Vs{D}(\Vs{u}^{n+1})) +\nabla p^{n+1} = {\Vs{f}}^{n+1},   \\
&\nabla \cdot \Vs{u}^{n+1} = 0
\end{aligned}
\right.
\quad \mbox{ on } \Omega,
\end{equation}
% \begin{equation}\label{Prob_TimeDiscretization}
% \left \{
% \begin{array}{rcc}
% \displaystyle\rho\frac{\Vs{u}^{n+1} - \Vs{u}^{n}}{\Delta t} + \rho\left(\Vs{u}^{n+1} \cdot \nabla \right)\Vs{u}^{n+1} - \nabla \cdot(2\nu^{n+1} \VD(\Vs{u}^{n+1})) + \nabla \, p^{n+1} &=&
% \Vs{f}^{n+1}\\ % \quad &\mbox{ in } \Omega, \\
% \nabla \cdot \Vs{u}^{n+1} &=& 0 \\ % &\mbox{ in } \Omega,\\
% \end{array}
% \right.
% \end{equation}
completed with the same initial and boundary conditions~\eqref{BC-Initial}-\eqref{BC-Limit}. Observe that since the viscosity could depend on the velocity (or pressure), \eqref{Prob_TimeDiscretization} contains possibly two non linearities. As for the first time step $t^1$, a simple backward Euler time step could be used. 

%%%%%%%%%%%%
% We must emphasizes that for most of these strategies (implicit, semi-explicit and explicit),
% the system \eqref{Prob_TimeDiscretization} will still be a coupled system. The velocity of the
% fluid $\Vs{u}$ depends on the temperature $\theta$ through the viscosity and the right member and 
% for the temperature, we have a convective term depending on the fluid velocity.
% Therefore an iterative scheme must be introduced to solve \eqref{Prob_TimeDiscretization} at each time step.

% For the rest of this section, we will present an existence result. In the next section, the determination
% of a solution of \eqref{Prob_TimeDiscretization} at each time step will be treated.
%The totally implicit approach has been retained since in all cases, a fixed-point will be needed (to deal with
%the coupling of the unknown). Furthermore the theoretical results presented in both sections are
%easily expanded to the semi-explicit (and even explicit) approach (corollaries \ref{coro_exitence} and
%\ref{coro_convergence}).
%%%%%%%%%%%%%%%

%\section{Fractional time-stepping and the projection scheme}\label{Projection} %\hspace{20cm}
\section{Toward a projection scheme for non--Newtonian fluid}

In this section, we follow the approach used in \cite{DetYak2017} leading to the creation of a more precise projection scheme. This new scheme is a generalization of the \textit{incremental projection}, as presented in \cite{GueMinShe2006} and its construction fellow the idea of the \textit{rotational projection} in \cite{TimMinVan1996}.

%Recall that, at time $t^{n+1}$, the incremental projection correspond to the 
Starting with the incremental projection, at time $t^{n+1}$, the method consists in the
following series of steps: solving the viscous nonlinear system (or Burgers' equation) \eqref{prediction1} gives a \textit{velocity prediction}, next is a projection step \eqref{correction_vitesse_pression_rotationnel} producing a \textit{pressure correction} and finally the pressure updating step \eqref{correction_pression}. 

\begin{equation}\label{prediction1}
\left \{
\begin{aligned}
\rho &D_t^\star\Vs{\tilde{\Vs{u}}} + \rho \left(\Vs{\tilde{\Vs{u}}} \cdot \nabla \right)\Vs{\tilde{\Vs{u}}} 
- \nabla \cdot(2\nu^{n+1}(\Vs{\tilde{\Vs{u}}}) \Vs{D}(\Vs{\tilde{\Vs{u}}})) +\nabla p^n = {\Vs{f}}^{n+1}
%\qquad \mbox{ on } \Omega,
\vspace*{4pt}\\
&\Vs{\tilde{\Vs{u}}} = 0
\qquad \mbox{ on } \Gamma_D
\end{aligned}
\right.
\end{equation}
\begin{equation}\label{correction_vitesse_pression_rotationnel}
\left \{
\begin{aligned}
\Vs{u}^{n+1} &= \Vs{\tilde{\Vs{u}}}-\displaystyle\frac{2\Delta t}{3\rho}\nabla \varphi \\
\nabla\cdot\Vs{u}^{n+1} &= 0  \\
\Vs{u}^{n+1}\cdot\Vs{n} &= 0 \qquad   \mbox{ on } \Gamma_D
\end{aligned}
\right.
\end{equation}
\begin{equation}\label{correction_pression}
p^{n+1} = p^n + \varphi
\end{equation}
with
$$
D_t^\star\Vs{\tilde{\Vs{u}}} = \displaystyle\frac{3\Vs{\tilde{\Vs{u}}}-4\Vs{u}^n+\Vs{u}^{n-1}}{2\Delta t} = D_t\Vs{u}^{n+1} - \frac{1}{\rho}\nabla\varphi.
$$
\begin{remark}\label{rmkInc}
We could choose to replace $\nu^{n+1}(\Vs{\tilde{\Vs{u}}})= \nu(t^{n+1},\Vs{\tilde{\Vs{u}}})$ by $\nu^{n+1}(\Vs{u}^n)$ in \eqref{prediction1}, making the non linearity related to the viscosity disappear. Obviously using $\nu^{n+1}(\Vs{u}^{n+1})$ makes the approach almost unusable as it would impose to solve \eqref{prediction1}--\eqref{correction_vitesse_pression_rotationnel} as a coupled system. 
\end{remark}

The first objection regarding \eqref{prediction1}--\eqref{correction_pression} is the fact that an artificial homogeneous Neumann boundary condition is enforced on the pressure. However, as illustrated in Figure~\ref{fig_acc_time_velo}--\ref{fig_acc_time_pre}, this inconsistent boundary condition as little effect on the overall precision of the scheme (see Rannacher et al.~\cite{Ran1992} for a detailed review of this point) and the $H^1$--norm of the error on the velocity prediction has a good behaviour. 

More importantly this splitting neglect information easily at our disposal (related to the shear rate of the velocity prediction). The approximation resulting from \eqref{prediction1}--\eqref{correction_pression} can be improved simply by recovering this information.
% propose an approximation of the pressure (and indirectly the velocity) that leaves behind relevant information based on the shear rate of $\Vs{\tilde{\Vs{u}}}$ and $\Vs{u}$. 
%
In \cite{TimMinVan1996} the rotational pressure correction scheme is proposed. This scheme gives a \textit{quasi-consistent} boundary condition for the pressure and better precision for the velocity and pressure, see \cite{GueMinShe2006}. However, this approach is only valid for \textbf{homogeneous viscosity}. 

The more general \textit{shear rate projection} proposed in \cite{DetYak2017},  makes it possible to achieve this improved accuracy for heterogeneous viscosity. But as presented it does not take into account the explicit dependency of the viscosity upon velocity, leaving again valuable information that could enrich the pressure approximation. 
Assuming that the viscosity is heterogeneous and depends on the velocity field $\Vs{u}$,  
we consider the new projection
\begin{equation}\label{prediction2}
\left \{
\begin{aligned}
\rho &D_t^\star\Vs{\tilde{\Vs{u}}}   + \rho \left(\Vs{\tilde{\Vs{u}}} \cdot \nabla \right)\Vs{\tilde{\Vs{u}}} 
- \nabla \cdot(2\nu^{n+1}(\Vs{\tilde{\Vs{u}}}) \Vs{D}(\Vs{\tilde{\Vs{u}}})) +\nabla p^n = {\Vs{f}}^{n+1}  
%\qquad \mbox{ on } \Omega,
\\
&\Vs{\tilde{\Vs{u}}} = 0\qquad \mbox{ on } \Gamma_D
\end{aligned}
\right.
\end{equation}
\begin{equation}\label{correction_vitesse_pression_rotationnel2}
\left \{
\begin{aligned}
\Vs{u}^{n+1} &= \Vs{\tilde{\Vs{u}}}-\displaystyle\frac{2\Delta t}{3\rho}\nabla \varphi \\
\nabla\cdot\Vs{u}^{n+1} &= 0  \\
\Vs{u}^{n+1}\cdot\Vs{n} &= 0 \qquad   \mbox{ on }\Gamma_D
\end{aligned}
\right.
\end{equation}
\begin{equation}\label{maj_pression}
p^{n+1} = p^n + \varphi + \psi
\end{equation}
% \begin{equation}\label{prediction2}
% \left \{
% \begin{array}{rcll}
% \displaystyle \rho \left(\frac{{3\Vs{\tilde{\Vs{u}}}}^{n+1}-4\Vs{u}^n+2\Vs{u}^{n-1}}{2\Delta t}\right)   + \rho \left(\Vs{\tilde{\Vs{u}}}^{n+1} \cdot \nabla \right)\Vs{\tilde{\Vs{u}}}^{n+1}\qquad& & &\\
% - \nabla \cdot(2\nu(\Vs{\tilde{\Vs{u}}}^{n+1}) \VD(\Vs{\tilde{\Vs{u}}}^{n+1})) +\nabla p^n &=& {\Vs{f}}^{n+1}&  \, \mbox{ on } \Omega,\\
% \Vs{\tilde{\Vs{u}}}^{n+1} &=& 0& \, \mbox{ on } \Gamma_D
% \end{array}
% \right.
% \end{equation}
% \begin{equation}\label{correction_vitesse_pression_rotationnel2}
% \left \{
% \begin{array}{rcll}
% \Vs{u}^{n+1} &=& \Vs{\tilde{\Vs{u}}}^{n+1}-\displaystyle\frac{2\Delta t}{3\rho}\nabla \varphi^{n+1}& \\
% p^{n+1} &=& p^n + \varphi^{n+1} - \psi^{n+1}\\
% \Vs{u}^{n+1} &=& 0 &\,   \mbox{ on } \Gamma_D
% \end{array}
% \right.
% \end{equation}
We will now establish the equation characterizing $\psi$ to obtain a splitting. From \eqref{correction_vitesse_pression_rotationnel2} and \eqref{maj_pression}
%, using the properties of the rotational operator 
we have
%$$\nabla\times \Vs{\tilde{\Vs{u}}}^{n+1} = \nabla\times\Vs{u}^{n+1}$$
$$ \nabla p^n = \nabla p^{n+1}- 3\rho\left(\frac{\Vs{\tilde{\Vs{u}}} - \Vs{u}^{n+1}}{2\Delta t}\right) 
- \nabla\psi.$$
Replacing in \eqref{prediction2} we get
$$
\begin{aligned}
\displaystyle \rho D_t\Vs{u}^{n+1} &  + \rho \left(\Vs{\tilde{\Vs{u}}} \cdot \nabla \right)\Vs{\tilde{\Vs{u}}} 
 - \nabla \cdot(2\nu^{n+1}(\Vs{\tilde{\Vs{u}}}) \Vs{D}(\Vs{\tilde{\Vs{u}}})) +\nabla p^{n+1} - \nabla\psi= {\Vs{f}}^{n+1}
\end{aligned}
$$
then we want $\psi$ satisfying
% $$
% \nabla\psi= \nabla \cdot(2\nu^{n+1} \VD(\Vs{\tilde{\Vs{u}}}^{n+1})) - \nabla \cdot(2\nu^{n+1} \VD(\Vs{u}^{n+1}))
% $$
% We will now establish the equation characterizing $\psi$. From \eqref{correction_vitesse_pression_rotationnel2}, using the properties of the ratational operator we have
% $$\nabla\times \Vs{\tilde{\Vs{u}}}^{n+1} = \nabla\times\Vs{u}^{n+1}$$
% $$ \nabla p^n = \nabla p^{n+1}- 3\rho\left(\frac{\Vs{\tilde{\Vs{u}}}^{n+1} - \Vs{u}^{n+1}}{2\Delta t}\right) 
% + \nabla\psi^{n+1}$$
% Replacing in \eqref{prediction2} we get
% $$
% \begin{aligned}
% \displaystyle \rho \left(\frac{{3\Vs{u}}^{n+1}-4\Vs{u}^n+2\Vs{u}^{n-1}}{2\Delta t}\right) &  + \rho \left(\Vs{\tilde{\Vs{u}}}^{n+1} \cdot \nabla \right)\Vs{\tilde{\Vs{u}}}^{n+1} \\
% & - \nabla \cdot(2\nu(\Vs{\tilde{\Vs{u}}}^{n+1}) \VD(\Vs{\tilde{\Vs{u}}}^{n+1})) +\nabla p^{n+1} + \nabla\psi^{n+1}= {\Vs{f}}^{n+1}
% \end{aligned}
% $$
% for $\psi$ satisfying
\begin{equation}\label{eqPsi0}
\nabla\psi= \nabla \cdot(2\nu^{n+1}(\Vs{u}^{n+1}) \Vs{D}(\Vs{u}^{n+1}))- \nabla \cdot(2\nu^{n+1}(\Vs{\tilde{\Vs{u}}}) \Vs{D}(\Vs{\tilde{\Vs{u}}})).
\end{equation}
% Theorem 2.9 in \cite{GirRav1979} insure the existence (and uniqueness up to a constant) of such $\psi$. Using the Agmon-Douglis-Nirenberg theorem (see \cite{BoyFab2006}), \eqref{eqPsi} is well defined and we can establish that $\psi$ is atleast as regular as $p$.
% %\end{remark}
%
% In \eqref{eqPsi} the right handside is known, and this splitting has two Poisson problems to solve. 
Taking the divergence, we have
% $$
% \Delta\psi^{n+1}= \nabla\cdot\nabla \cdot\left(2\nu(\Vs{\tilde{\Vs{u}}}^{n+1}) \VD(\Vs{\tilde{\Vs{u}}}^{n+1}) - 2\nu(\Vs{u}^{n+1})\VD(\Vs{u}^{n+1})\right)$$
\begin{equation}\label{eqPsi}
\Delta\psi= \nabla\cdot\nabla \cdot\left(2(\nu^{n+1}(\Vs{u}^{n+1}) - \nu^{n+1}(\Vs{\tilde{\Vs{u}}}))\Vs{D}(\Vs{\tilde{\Vs{u}}})\right) -\nabla\cdot\nabla \cdot\left( 2\nu^{n+1}(\Vs{u}^{n+1})\Vs{D}(\Vs{\tilde{\Vs{u}}} - \Vs{u}^{n+1})\right).
\end{equation}

%%%%%%%%%%%%%%%
\noindent The right hand side of this equation is known, as for the boundary conditions for $\psi$, the compatibility condition for the Poisson equation gives us a Neumann condition (see \eqref{EqPsi}). For $\Vs{f}$ sufficiently regular, using the Agmon-Douglis-Nirenberg theorem (see \cite{BoyFab2006}), the solution of this Poisson problem, is well defined (see Remark~ \ref{remarkPoisson}) and we can establish that $\psi$ is of the same nature and as regular as the pressure.

Since the convective term is on $\Vs{\tilde{\Vs{u}}}$, as in \cite{TimMinVan1996, GueMinShe2006, DetYak2017}, from \eqref{prediction2} and \eqref{eqPsi0} we get a \textit{nearly} consistent boundary condition for $p^{n+1}$ on $\Gamma_D$. A consistent condition would be obtained provided the non linear convective term in the momentum equation is neglected; which correspond to an unsteady Stokes problem.
$$
\nabla p^{n+1}\cdot \Vs{n} = \left({\Vs{f}}^{n+1} - \rho D_t\Vs{u}^{n+1}   + \rho \left(\Vs{\tilde{\Vs{u}}} \cdot \nabla \right)\Vs{\tilde{\Vs{u}}}
- \nabla \cdot(2\nu^{n+1}(\Vs{u}^{n+1}) \Vs{D}(\Vs{u}^{n+1}))\right)\cdot\Vs{n}.
$$
%%%%%%%%%%%%%%

% As for the boundary conditions for $\psi$, the compatibility condition for the Poisson equation gives us a Neumann condition along $\Gamma_D$. 

\begin{remark}
% First some algebraic identity related to a field $\V\in \mathbb R^d,~ d=$  2 or 3 and a scalar function $\psi$
% \begin{equation}\label{rotrot}
% \nabla \times \nabla \times \V = \nabla(\nabla\cdot\V) - \Delta \V
% \end{equation}
% \begin{equation}\label{divDu}
% \nabla\cdot\left(2\VD(\V))\right) = \Delta\V + \nabla(\nabla\cdot \V)
% \end{equation}
% \begin{equation}\label{divnu}
% \begin{array}{rcll}
% \nabla\cdot\left(2\psi\VD(\V))\right) &=& 2\nabla\psi \VD(\V) + \psi \Delta\V + \psi \nabla(\nabla\cdot \V)\\
% &=& 2\nabla\psi \VD(\V) - \psi (\nabla \times \nabla \times \V)  + 2\psi\nabla(\nabla\cdot \V)
% \end{array}
% \end{equation}
Obviously in case of solenoidal vector fields the right hand side of \eqref{eqPsi} could be further simplified. However, a priori, the velocity prediction $\Vs{\tilde{\Vs{u}}}$ is not divergence free.
From \eqref{eqPsi}, $\psi$ is composed of two corrections related to the difference between $\Vs{\tilde{\Vs{u}}}$ and $\Vs{u}^{n+1}$: the first is associated with the variation of viscosity, the second one related to the variation of shear rate.
\end{remark}
%
%
% Using the Agmon-Douglis-Nirenberg theorem (see \cite{BoyFab2006}), \eqref{eqPsi} is well defined and we can establish that $\psi$ is atleast as regular as $p$.
%\end{remark}
% In \eqref{eqPsi} the left handside is known, as for the boundary conditions for $\psi$, the compatibility condition for the Poisson equation gives us a Neumann condition along $\Gamma_D$. 
% Finally, going back to \eqref{correction_vitesse_pression_rotationnel2}, we get on $\Gamma_D$, as in \cite{XXX, YYY,DetYak2017}, a quasi-consistent condition on $p^{n+1}$.
% $$
% \nabla p^{n+1}\cdot \Vs{n} = \left({\Vs{f}}^{n+1} - \rho D_t\Vs{u}^{n+1}   + \rho \left(\Vs{\tilde{\Vs{u}}} \cdot \nabla \right)\Vs{\tilde{\Vs{u}}}
% - \nabla \cdot(2\nu^{n+1}(\Vs{u}^{n+1}) \VD(\Vs{u}^{n+1}))\right)\cdot\Vs{n}.
% $$

%
%
%
\section{Strategies to solve the shear rate projection for non--Newtonian fluids}\label{Strategies}

The prediction step \eqref{prediction2} imply solving a problem containing two  nonlinearities: the convective term and the viscous term (since $\nu$ depends on the velocity and possibly the pressure, etc.). %Various strategies can be considered. 
Introducing $\Vs{\tilde{\Vs{u}}}^\star$ and $\nu^\star$ to lighten the notation, 
%les termes dans 
\eqref{prediction2} is written
%$$(\Vs{\tilde{\Vs{u}}}^\star\cdot\nabla)\Vs{\tilde{\Vs{u}}}\qquad\text{and}\qquad \nabla\cdot(2\nu^\star\VD(\Vs{\tilde{\Vs{u}}}))$$
\begin{equation}\label{prediction3}
\left \{
\begin{aligned}
&3\rho \frac{ \Vs{\tilde{\Vs{u}}}}{2\Delta t} +\rho\left( \Vs{\tilde{\Vs{u}}}^\star\cdot \nabla \right) \Vs{\tilde{\Vs{u}}}
- \nabla \cdot (2\nu^\star \Vs{D}(\Vs{\tilde{\Vs{u}}}) ) 
= \Vs{f}^{n+1} - \nabla p^n  +\rho\frac{4\Vs{u}^{n}-\Vs{u}^{n-1}}{2\Delta t}
\\
&\tilde{\Vs{u}} = \textbf{0} \quad \mbox{on}\quad \Gamma_D %\partial \Omega. 
\end{aligned} 
\right.
\end{equation}
which encompass various strategies to solve \eqref{prediction2} according to the definition given for $\Vs{\tilde{\Vs{u}}}^\star$ and $\nu^\star$.

Defining $\Vs{\tilde{\Vs{u}}}^\star$ and $\nu^\star$ as expression not depending on  $\Vs{u}^{n+1}$ would lead to "explicit" strategies, conversely using expression depending on $\Vs{u}^{n+1}$ or $\Vs{\tilde{\Vs{u}}}$ would be called "implicit" strategies. These implicit approaches correspond to some fixed point methods applied to solve \eqref{prediction3}: at time step $t^{n+1}$, given $\Vs{\tilde{\Vs{u}}}_0$ an initial approximation of $\Vs{\tilde{\Vs{u}}}$ solution of \eqref{prediction2},  until convergence, we seek $\Vs{\tilde{\Vs{u}}}_{k+1}$ solution of
\begin{equation}\label{linconv}
\left \{
\begin{aligned}
&3\rho \frac{ \Vs{\tilde{\Vs{u}}}_{k+1}}{2\Delta t} +\rho\left( \Vs{\tilde{\Vs{u}}}^\star\cdot \nabla \right) \Vs{\tilde{\Vs{u}}}_{k+1}
- \nabla \cdot (2\nu^\star \Vs{D}(\Vs{\tilde{\Vs{u}}}_{k+1}) ) 
= \Vs{f}^{n+1} - \nabla p^n  +\rho\frac{4\Vs{u}^{n}-\Vs{u}^{n-1}}{2\Delta t}
\\
&\tilde{\Vs{u}}_{k+1} = \textbf{0} \qquad \mbox{on}\quad \Gamma_D. %\partial \Omega.
\end{aligned} 
\right.
\qquad k=1,2,\dots
\end{equation}

If the choice of treatment can impact the overall order of convergence of the time discretization, it has no fundamental consequence on the viscous correction presented. We propose a few examples of treatment of both terms.

\subsection{Treatment of convection term}
Explicit approaches consists in using extrapolation (backward in time) for $\Vs{\tilde{\Vs{u}}}^\star$,
leading %. These approaches leads 
to a complete linearization of the convective term. The simplest approach consist in using, at time step $t^{n+1}$, %the extrapolation
$$\Vs{\tilde{\Vs{u}}}^* = \Vs{u}^n.$$
A more effective and accurate strategy is the use of the Richardson extrapolation
$$\Vs{\tilde{\Vs{u}}}^* = 2\Vs{u}^n - \Vs{u}^{n-1}.$$
For implicit approaches an extrapolation (backward with respect to the iterative steps) of the convective term is used. 
The choice of $\Vs{\tilde{\Vs{u}}}^\star$ will alter the behaviour of the iterative process. The simplest extrapolation consists, at time step $t^{n+1}$ for the iterative step $k+1$, in using 
$$\Vs{\tilde{\Vs{u}}}^\star = \Vs{\tilde{\Vs{u}}}_k.$$
Once again a Richardson extrapolation produce another strategy 
%(usually more efficient then the basic Picard iteration)
$$\Vs{\tilde{\Vs{u}}}^\star = 2\Vs{\tilde{\Vs{u}}}_k - \Vs{\tilde{\Vs{u}}}_{k-1}.$$
Finally, %although it is a fixed point approach, 
Newton's method is treated separately as it would add terms to \eqref{linconv}. In this case the use of $\Vs{\tilde{\Vs{u}}}^\star$ can be dropped and~\eqref{linconv} is replaced by~\eqref{newtonconv}: at time step $t^{n+1}$, given $\Vs{\tilde{\Vs{u}}}_0$ an initial approximation of $\Vs{\tilde{\Vs{u}}}$ solution of \eqref{prediction2} and until convergence, we seek $\Vs{\tilde{\Vs{u}}}_{k+1}$ solution of
\begin{equation}\label{newtonconv}
\left \{
\begin{array}{l}
\begin{aligned}
3\rho \frac{ \Vs{\tilde{\Vs{u}}}_{k+1}}{2\Delta t} 
+\rho\left(
\left( \Vs{\tilde{\Vs{u}}}_k\cdot \nabla \right) \Vs{\tilde{\Vs{u}}}_{k+1}
+\left( \Vs{\tilde{\Vs{u}}}_{k+1}\cdot \nabla \right) \Vs{\tilde{\Vs{u}}}_{k}\right)
- \nabla \cdot (&2\nu^\star \Vs{D}(\Vs{\tilde{\Vs{u}}}_{k+1}) ) \\ 
= \Vs{f}^{n+1} - \nabla p^n  +\rho&\frac{4\Vs{u}^{n}-\Vs{u}^{n-1}}{2\Delta t}
+\rho\left( \Vs{\tilde{\Vs{u}}}_{k}\cdot \nabla \right) \Vs{\tilde{\Vs{u}}}_{k}
\end{aligned} \vspace*{4pt}\\
\tilde{\Vs{u}}_{k+1} = \textbf{0} \quad \mbox{on}\quad \Gamma_D. %\partial \Omega.
\end{array}
\right.
\qquad k=1,2,\dots
\end{equation}
\subsection{Treatment of the viscous term}
The same basic strategies can be applied to the non linear viscous term. Explicit approaches consists in using extrapolation (backward in time) for $\nu^\star$. The simplest approach 
%consist in using, 
at time step $t^{n+1}$, is %the extrapolation
$$\nu^\star = \nu(t^{n+1}, \Vs{u}^n)=\nu^{n+1}(\Vs{u}^n).$$
Two kind of Richardson extrapolations can be considered
$$\nu^\star = \nu(t^{n+1}, 2\Vs{u}^n - \Vs{u}^{n-1})\qquad\text{ or }\qquad \nu^\star = 2\nu(t^{n+1}, \Vs{u}^n) - \nu(t^{n+1},\Vs{u}^{n-1}).$$
For implicit approaches, backward extrapolations (with respect to the iterative steps) of the viscosity are introduced. 
The choice of $\nu^\star$ will alter the behaviour of the iterative process. The simplest extrapolation at time step $t^{n+1}$ for the iterative step $k+1$, is %in using 
$$\nu^\star = \nu(t^{n+1}, \Vs{\tilde{\Vs{u}}}_k)=\nu^{n+1}(\Vs{\tilde{\Vs{u}}}_k).$$
Once again a Richardson extrapolation produce other strategies 
$$\nu^\star = \nu^{n+1}(2\Vs{\tilde{\Vs{u}}}_k - \Vs{\tilde{\Vs{u}}}_{k-1})\qquad\text{ or }\qquad \nu^\star = 2\nu(t^{n+1}, \Vs{\tilde{\Vs{u}}}_k) - \nu(t^{n+1},\Vs{\tilde{\Vs{u}}}_{k-1}).$$
Finally, assuming we have sufficient regularity of $\nu$ with respect to $\Vs{u}$, we can apply Newton's method. Introducing 
$$\delta_\nu^{n+1}(\Vs{u}) = \nabla_\Vs{u} \nu(t^{n+1},\Vs{u})$$ 
the gradient of $\nu$ with respect to $\Vs{u}$ at time $t^{n+1}$. Dropping the $\nu^\star$ notation  and replacing~\eqref{linconv} by~\eqref{newtonvisc} we get: given $\Vs{\tilde{\Vs{u}}}_0$ an initial approximation of $\Vs{\tilde{\Vs{u}}}$ solution of \eqref{prediction2}, we seek $\Vs{\tilde{\Vs{u}}}_{k+1}$ solution of
\begin{equation}\label{newtonvisc}
\left \{
\begin{array}{l}
\begin{aligned}
3\rho \frac{ \Vs{\tilde{\Vs{u}}}_{k+1}}{2\Delta t} 
+\rho \left( \Vs{\tilde{\Vs{u}}}^\star\cdot \nabla \right) \Vs{\tilde{\Vs{u}}}_{k+1}
- \nabla \cdot (&2\nu^{n+1}(\Vs{\tilde{\Vs{u}}}_k) \Vs{D}(\Vs{\tilde{\Vs{u}}}_{k+1}) + 2(\delta_\nu^{n+1}(\Vs{\tilde{\Vs{u}}}_k) \cdot\Vs{\tilde{\Vs{u}}}_{k+1})\Vs{D}(\Vs{\tilde{\Vs{u}}}_{k})) \\ 
&= \Vs{f}^{n+1} - \nabla p^n  +\rho\frac{4\Vs{u}^{n}-\Vs{u}^{n-1}}{2\Delta t}
-\nabla\cdot(2(\delta_\nu^{n+1}(\Vs{\tilde{\Vs{u}}}_k)\cdot\Vs{\tilde{\Vs{u}}}_k)\Vs{D}(\Vs{\tilde{\Vs{u}}}_{k}))
\end{aligned} \vspace*{4pt}\\
\tilde{\Vs{u}}_{k+1} = \textbf{0} \quad \mbox{on}\quad \partial \Omega. 
\end{array}
\right.
\end{equation}

Both non linear terms can be treated explicitly, implicitly or by "hybrid" strategy (using an explicit approach for one term and an implicit for the other). Notice that the use of explicit strategies for both terms makes \eqref{prediction3} a linear equation, in all other cases we have a fixed point loop of general form \eqref{linconv} to consider (including a possible combined expression composed of \eqref{newtonconv} and \eqref{newtonvisc}). 
\begin{remark}
It is tempting to follow the same idea for the pressure, replacing $p^n$ in~\eqref{prediction2} and~\eqref{maj_pression}  by $p^\star = 2p^n-p^{n-1}$ a second order extrapolation. However, as observed in~\cite{GueMinShe2006}, in that case the stability seems to necessitate a lower bound on $\Delta t$. 
\end{remark}
%%%%%%%%%%%%%%%%%%%%%%%%
\subsection{A shear rate projection algorithm for heterogeneous viscosity}

%We propose a procedure based on 
% Various approach can be used for the nonlinear convective term in \eqref{prediction2}. 
%This choice will influence the overall speed of the procedure but also its accuracy. 
%Step 2 correspond to the use of a simple
%COR:Using a simple %
We now summarize these different strategies in an algorithm, keeping in mind that certain choice for $\nu^\star$ and $\tilde{\Vs{u}}^\star$ can modify the expression~\eqref{AlgNSDivu} in Step 2 below, or make it a linear system (in case of explicit expression).
As a generic strategy for solving \eqref{prediction2}  we get the following algorithm: at each time step $t^{n+1}$, 

\begin{enumerate}
\item {{\it Initialization:}} 
$$
\Vs{\tilde{\Vs{u}}}_0=\Vs{u}^n.
$$
\item {\it{Solving the non linear equation in $\Vs{\tilde{\Vs{u}}}^{n+1}$}:}

\qquad until convergence, compute $\Vs{\tilde{\Vs{u}}}_{k+1}$ solution of:
\begin{equation}\label{AlgNSDivu}
\left \{
%\begin{array}{l}
\begin{aligned}
3&\rho \frac{ \Vs{\tilde{\Vs{u}}}_{k+1}}{2\Delta t} +\rho\left( \Vs{\tilde{\Vs{u}}}^\star\cdot \nabla \right) \Vs{\tilde{\Vs{u}}}_{k+1}
- \nabla \cdot ( 2\nu^\star \Vs{D}(\Vs{\tilde{\Vs{u}}}_{k+1}) ) 
= \Vs{f}^{n+1} - \nabla p^n  +\rho\frac{4\Vs{u}^{n}-\Vs{u}^{n-1}}{2\Delta t}
\vspace*{4pt}\\
&\Vs{\tilde{\Vs{u}}}_{k+1} = \textbf{0} \quad \mbox{on}\quad \Gamma_D. 
\end{aligned} 
%\end{array}
\right.
\end{equation}
\qquad When converged put $\Vs{\tilde{\Vs{u}}}^{n+1} = \Vs{\tilde{\Vs{u}}}_{k+1}$.
\item {{\it Projection step :}}  

\qquad Compute  $\varphi$ solution of Poisson problem:
%COR:Compute  $\varphi$ solution of Poisson problem with suitable boundary conditions:
\begin{equation}\label{EqPhi}
\Delta \varphi = \nabla \cdot \left(\frac{3\rho}{2\Delta t} \Vs{\tilde{\Vs{u}}}^{n+1}\right)
\qquad\qquad \varphi = 0\quad\text{on }\Gamma_N %=\partial\Omega\backslash\Gamma_D
\end{equation}
\item {{\it  Velocity correction :}}
$$\Vs{u}^{n+1} = \Vs{\tilde{\Vs{u}}}^{n+1}- \frac{2\Delta t}{3\rho} \nabla \varphi$$%\qquad\qquad \Vs{u}^{n+1} = 0\quad\text{on }\Gamma_D$$
\item {{\it Shear rate projection step :}}  

\qquad Compute  $\psi$ solution of Poisson problem:
\begin{equation}\label{EqPsi}
\left \{
\begin{aligned}%{array}{lcl}
&\Delta\psi = 
%\nabla \cdot((2\nu^{n+1}(\Vs{u}^{n+1}) \Vs{D}(\Vs{u}^{n+1})) - (2\nu^{n+1}(\Vs{\tilde{\Vs{u}}}) \Vs{D}(\Vs{\tilde{\Vs{u}}})))
\nabla\cdot\nabla \cdot\left( 
2\nu^{n+1}(\Vs{u}^{n+1})\Vs{D}(\Vs{u}^{n+1}) - 2\nu^\star\Vs{D}(\Vs{\tilde{\Vs{u}}}^{n+1}) 
\right) \\[0.5em]
& \nabla\psi\cdot\Vs{n} = \left(\nabla \cdot\left(
2\nu^{n+1}(\Vs{u}^{n+1})\Vs{D}(\Vs{u}^{n+1}) - 2\nu^\star\Vs{D}(\Vs{\tilde{\Vs{u}}}^{n+1})
\right)\right)\Vs{n}\qquad\text{on }\partial\Omega
\end{aligned}%{array}
\right.
\end{equation}
\item {{\it Pressure correction :}}  
$$p^{n+1} = p^n + \varphi+\psi.$$
\end{enumerate}

\begin{remark}\label{PremierPas}
For the first time step ($n=1$) a simple backward Euler or Crank-Nicolson scheme is used.
%, this slightly modifies the expressions in Step 2, 3 and 4 at the first time step. 
Step 4 and 6 have to be interpreted as $L^2$-projections in the velocity and pressure space respectively. 
\end{remark}
\begin{remark} \label{remarkPoisson}%\ \\
%\begin{enumerate}[i)]%label=(\roman*)]
%\item 
For $\Gamma_N=\emptyset$ the Poisson problem in step 3 is ill-posed and we seek a solution in $H^1(\Omega)\backslash\mathbb{R}$.
For $\Gamma_N\neq\emptyset$ the homogeneous Dirichlet boundary condition in \eqref{EqPhi} is imposed by the Helmoltz decomposition (as in \cite{GueMinShe2005}). 
From \cite{GirRav1986}, step 5 has 
%the divergence formula is used to limit the order of the derivatives and 
a unique solution in $H^1(\Omega)\backslash\mathbb{R}$ and its variational form is 
$$(\nabla\psi,\nabla v) = \left(\nabla\cdot\left(
2\nu^{n+1}(\Vs{u}^{n+1})\Vs{D}(\Vs{u}^{n+1}) - 2\nu^\star\Vs{D}(\Vs{\tilde{\Vs{u}}}^{n+1})
\right), \nabla v\right)\qquad\qquad\forall v \in H^1(\Omega).$$
%$$\int\limits_{\Omega}\nabla\psi\cdot\nabla v\ d\Omega = \int\limits_{\Omega}\nabla\cdot(2\nu \mathbf{D}(\mathbf{u}-\mathbf{\tilde u})))\cdot \nabla v\ d\Omega\qquad\qquad\forall v \in H^1(\Omega)$$
%The equality for the velocity and pressure correction 
%
%\item 
From a computational point of view, \eqref{EqPhi} and \eqref{EqPsi} can be viewed as solving the same problem twice with two different right hand side. For methods such as the finite element, strategies reducing the computational effort involved in solving both solutions are available. 
%
%\end{enumerate}
\end{remark}
\begin{remark}
The basic incremental projection \cite{God1979}, consists in ignoring step 5 and putting $\psi \equiv 0$ which gives the original system \eqref{prediction1}-\eqref{correction_vitesse_pression_rotationnel}. 
²For an homogeneous viscosity, solving \eqref{EqPsi} is not required as in this case, thanks to the properties of the rotational operator, we have
$$
\nabla\cdot\nabla\psi = -\nabla\cdot(2\nu^{n+1} \nabla(\nabla\cdot\Vs{\tilde{\Vs{u}}}^{n+1}))\quad\Rightarrow \psi = - 2 \nu^{n+1} \nabla\cdot\Vs{\tilde{\Vs{u}}}^{n+1}
$$
which is the rotational projection \cite{TimMinVan1996}. In case of a heterogeneous viscosity, not depending on $\Vs{u}$, obviously the first term in \eqref{eqPsi} will cancel, leaving only the second term resulting in the projection presented in \cite{DetYak2017}.
Finally, in step 6 (pressure correction), applying a "relaxation" factor $\alpha \in ]0,1]$ to $\psi$ we get a generalized version of the "Chorin-Uzawa" scheme as described by Rannacher in \cite{Ran2000}.  
%\end{enumerate}
\end{remark}

\section{Finite element discretization}
For the spatial approximation, the finite element method is used. First, we introduce a mesh $\mathcal{T}_{h} = \{\top\}$ of simplicial $\top$ partitioning $\Omega$:
\begin{equation*}\label{mesh}
      \overline{\Omega}_h := \bigcup_{\top \in \mathcal{T}_{h}}\top\,\subseteq \overline{\Omega}
\end{equation*}
with the usual restrictions on $\mathcal{T}_{h}$ (see \cite{Bat1996} or \cite{CiaLun2009} for details). We emphasize that at no point in what follows are we making assumptions on the nature of the mesh. This partitioning of $\Omega$ could be composed of triangles, quadrangles or both  for $\Omega\in\mathbb{R}^2$ or tetrahedra, hexahedra or both for $\Omega\in\mathbb{R}^3$.  
% in such a way that the intersection between two tetrahedra, when non-empty, is assumed to be a vertex, a segment, or a triangle of both elements. The triangulation
% $\mathcal{T}_{h}$ is assumed regular, that is there exists $\rho \geq 0$ such that for any $\top \in \mathcal{T}_{h}$,
% \begin{eqnarray}
%       \frac{h_{\top}}{\rho_{\top}} \leq \rho,
% \end{eqnarray}
% where $\rho_{\top}$ denotes the radius of the inscribed ball in $\top$ and $h_{\top}$ is the diameter of the simplex $\top$. Let us denote by $h$ the mesh size, defined by
% \begin{eqnarray}\label{meshsize}
%       h := \max_{\top \in \mathcal{T}_{h}}h_{\top}.
% \end{eqnarray}
From the family $\mathcal{T}_{h}$ of partition of the domain $\Omega$ (indexed by $h$), we construct the family of finite dimensional vector spaces. 
\begin{eqnarray*}
    \V_{h} = \left\{\Vv_h \in (C^{0}\left(\overline{\Omega_h}\right))^d;\, {\Vv_h}_{|\top} \in P_{k_{\Vu}}\left(\top\right) , \, \forall \, \top \in  \mathcal{T}_{h}, \, \Vv_h = 0 \, \mbox{ on } \Gamma_{D} \right\}
\end{eqnarray*}
% \begin{eqnarray*}
%     \tilde V_{h} = \left\{\Vv_h \in \left(\mathscr{C}^{0}\left(\overline{\Omega_h}\right)\right)^d;\, {\Vv_h}_{|\top} \in \mathscr{P}_{k}\left(\top\right) , \, \forall \, \top \in  \mathcal{T}_{h}, \, v_h = 0 \, \mbox{ on } \Gamma_{D} \right\}
% \end{eqnarray*}
$$M_h= \left\{\vp_h \in C^{0}\left(\overline{\Omega_h}\right);\, {\vp_h}_{|\top} \in P_{k_p}\left(\top\right) , \, \forall \, \top \in  \mathcal{T}_{h}, \right\}$$
$$M^{\varphi}_h= \left\{\varphi_h \in C^{0}\left(\overline{\Omega_h}\right);\, {\varphi_h}_{|\top} \in P_{k_\varphi}\left(\top\right) , \, \forall \, \top \in  \mathcal{T}_{h},  \, \varphi_h = 0 \, \mbox{ on } \Gamma_{N}\right\}$$
$$M^{\psi}_h= \left\{\psi_h \in C^{0}\left(\overline{\Omega_h}\right);\, {\psi_h}_{|\top} \in P_{k_\psi}\left(\top\right) , \, \forall \, \top \in  \mathcal{T}_{h}\right\}$$
where $P_{k_\alpha}(\top)$ is the space of polynomials of degree less or equal to $k_\alpha$ on $\top$. The index "$h$" is used to denote the different spatial approximations. The approximation of $f^n(\Vs{x})$, a function $f(t,\Vs{x})$ at time $t^n$, will be denoted $\text{f}^n_h$. At each time step $t^{n+1}$, we are now seeking approximations 
$$\pVu_h^{n+1},\,\Vu_h^{n+1}\in \V_h,\ \vp_h^{n+1}\in M_h,\ \varphi_h\in M_h^{\varphi},\ \psi_h\in M_h^{\psi}.$$ 

Since the discrete inf--sup condition must be respected for $(\Vu_h,\vp_h)$ \cite{GirRav1986, GueMinShe2006,BofBreFor2013}, the prediction equation \eqref{AlgNSDivu} 
(recall that \eqref{newtonconv}, \eqref{newtonvisc} or a composite expression from both could replace \eqref{AlgNSDivu}) 
is discretized using specific choice of interpolations, typically $k_p = k_\Vu -1$, see~\cite{BofBreFor2013} for the details and other choices. 
Here a second order Taylor-Hood $P_2-P_1$ element was chosen for $\V_h$ and $M_h$. %(see~\cite{BofBreFor2013} for the details and other choices), therefore $k_\Vu = 2$ and $k_p = 1$.
%
% The prediction equation \eqref{AlgNSDivu} is discretized using a second order ($O (h^2)$) Taylor-Hood $P_2-P_1$ element (see Brezzi and Fortin~\cite{BreFor1991}). Since the discrete inf-sup condition must be respected for $(u,p)$ \cite{GirRav1986, GueMinShe2006}. 
Regarding the approximation of $\varphi$ and $\psi$, a priori, they could be of arbitrary interpolating degree. Even though steps 4 and 6 are interpreted as projections in $\V_h$ and $M_h$, we must take into account that $\varphi$ and $\psi$ are solutions of Poisson problems and that they act as velocity and pressure corrections. In this work, $\varphi_h$ and $\psi_h$ were chosen to be quadratic ($k_\varphi=k_\psi=2$).

Following the usual approach, the Galerkin formulation based on $\V_h$, $M_h$, $M_h^{\varphi}$ and $M_h^{\psi}$ is used to generate finite dimensional systems corresponding to~\eqref{AlgNSDivu}, \eqref{EqPhi} and \eqref{EqPsi}. As underlined in Remark~\ref{remarkPoisson}, the discrete version of Step~3 and Step~5 must be treated in a special manner (at least one discrete system is under-determined) and a unique solution is obtained through a supplementary computation in the last step (see Step~6$_h$). To complete the spatial discretization of the shear rate algorithm a proper interpretation of Step~4 and Step~6 must be given. Here we chose the most natural approach which consist in the use of a $L^2$--projection of $\varphi_h$ and $\psi_h$ in $M_h$ (Step~6) and of $\nabla\varphi_h$ in $\V_h$ (Step~4) (see Remark~\ref{projHier}). We are now in measure of giving a completely discrete generic version of the shear rate projection algorithm:

%%%%%%%%%%%%%%%%%
\begin{enumerate}[1$_h$.]
\item {{\it Initialization:}} 
$$
\pVu_0=\Vu^n_h.
$$
\item {\it{Solving the non linear prediction equation}:} 

\qquad until convergence, compute $\pVu_{k+1}\in \V_h$ \textit{solution of}
\begin{equation}\label{AlgNSDivu_h}
\begin{aligned}
\int\limits_{\Omega_h} 3\rho \frac{ \pVu_{k+1}}{2\Delta t}\cdot\Vv\, d\Omega 
+\int\limits_{\Omega_h} \rho\left( \pVu^\star\cdot \nabla \right) \pVu_{k+1}\cdot\Vv\, d\Omega
&+ \int\limits_{\Omega_h} ( 2\nu^\star \Vs{D}(\pVu_{k+1}) )\cdot \nabla\Vv\, d\Omega\\  
= \int\limits_{\Omega_h} \Vs{f}^{n+1}\Vv\, d\Omega
- \int\limits_{\Omega_h} &\nabla \vp^n_h\cdot\Vv\, d\Omega  + \int\limits_{\Omega_h} \rho\frac{4\Vu^{n}_h-\Vu^{n-1}_h}{2\Delta t}\cdot\Vv\, d\Omega
\end{aligned} 
\qquad\forall \Vv\in \V_h
\end{equation}
\qquad When converged put $\pVu_h^{n+1} = \pVu_{k+1}$.\vspace*{0.5em}
\item {{\it Projection step :}}  

\qquad Compute  $\varphi_h\in M_h^{\varphi}$ solution of
\begin{equation}\label{EqPhi_h}
\int\limits_{\Omega_h} \nabla\varphi_h\cdot\nabla\zeta\, d\Omega = 
\int\limits_{\Omega_h} \frac{3\rho}{2\Delta t} \, \pVu_h^{n+1}\cdot\nabla \zeta\, d\Omega
\qquad\forall \zeta\in M_h^{\varphi}\end{equation}
\item {{\it  Velocity correction :}}

\qquad Compute  $\Vu^{n+1}_h\in \V_h$ solution of :
$$\int\limits_{\Omega_h} \Vu^{n+1}_h\Vv\, d\Omega = \int\limits_{\Omega_h} \pVu_h^{n+1}\cdot\Vv\, d\Omega - \int\limits_{\Omega_h} \frac{2\Delta t}{3\rho} \nabla \varphi_h\cdot\Vv\, d\Omega \qquad\forall \Vv\in \V_h$$
\item {{\it Shear rate projection step :}}  

\qquad Compute  $\psi_h\in M_h^{\psi}$ solution of
%%%%
% \begin{equation}\label{Projection_divPsi_h}
% \int\limits_{\Omega_h} \Vb{divS}_h\xi\, d\Omega = 
% \int\limits_{\Omega_h} \left( \nabla\cdot\Vb{S}_h\right) \xi\, d\Omega\qquad\qquad\forall \xi\in M_h
% \end{equation}
%%%%
\begin{equation}\label{EqPsi_h}
\int\limits_{\Omega_h} \nabla\psi_h\cdot\nabla \zeta\, d\Omega = 
\int\limits_{\Omega_h}
\left( \nabla\cdot\VS_h\right)
%\Vb{divS}_h
%\left(2\nu^\star\Vs{D}(\pVu_h) - 2\nu^{n+1}(\Vu^{n+1}_h)\Vs{D}(\Vu^{n+1}_h)\right)
\cdot\nabla \zeta\, d\Omega\qquad\qquad\forall \zeta\in M_h^{\psi}
\end{equation}
with
%\begin{equation}\label{Projection_Psi_h}
$\VS_h = 
\left(2\nu^\star\Vs{D}(\pVu_h^{n+1}) - 2\nu^{n+1}(\Vu^{n+1}_h)\Vs{D}(\Vu^{n+1}_h)\right)$
in $\left(M_h\right)^{d\times d}.$\vspace*{0.5em}
%\end{equation}
% \begin{equation}\label{Projection_Psi_h}
% \int\limits_{\Omega_h} \VS_h\VT\, d\Omega = 
% \int\limits_{\Omega_h} \left(2\nu^\star\Vs{D}(\pVu_h^{n+1}) - 2\nu^{n+1}(\Vu^{n+1}_h)\Vs{D}(\Vu^{n+1}_h)\right)\cdot\VT\, d\Omega\qquad\qquad\forall \VT\in \left(M_h\right)^{d\times d}
% \end{equation}
%%%%%%%%%%%%%%%%%%%%%%%%%%%%%%%%%%%%%%
\item {{\it Pressure correction :}}  \label{step6h}

\qquad Compute  $\vp^{n+1}_h\in M_h$ solution of 
$$\int\limits_{\Omega_h} \vp_h^{n+1}\zeta\, d\Omega = \int\limits_{\Omega_h} \vp_h^n\zeta\, d\Omega + \int\limits_{\Omega_h} \left(\varphi_h+\psi_h - \mu(\varphi_h,\psi_h)\right)\zeta\, d\Omega
\qquad\forall \zeta\in M_h$$
\qquad with
$$
\mu(\varphi_h,\psi_h) = \frac{1}{mes(\Omega_h)} \begin{dcases}
%\displaystyle\int\limits_{\Omega_h}\, d\Omega}
%\, 
\displaystyle
\int\limits_{\Omega_h} \left(\varphi_h+\psi_h\right)\, d\Omega&\, \text{ if }\Gamma_N = \emptyset \\
%\, 
\displaystyle
\int\limits_{\Omega_h} \psi_h\, d\Omega&\, \text{ if } \Gamma_N \neq \emptyset 
\end{dcases}
$$
\end{enumerate}
%%%%%%%%%%%%%%%%%
%
%%===========
%=======================
%
\begin{remark}For $n=1$, the use of a backward Euler or Crank-Nicolson scheme slightly modifies the expressions in Step 2$_h$, 3$_h$ and 4$_h$. In Step 5$_h$, $\VS_h$ is defined through an $L^2$ projection. In Step 6$_h$, %correspond to $L^2$-projection on $\V_h$ and $M_h$ respectively. As for 
$\mu(\cdot,\cdot)$ correspond to an average on $\Omega_h$ of $\psi_h$ or $\varphi_h+\psi_h$. 
This correction gives $\psi_h\in M_h^{\psi}\backslash\mathbb{R}$ (and $\varphi_h\in M_h^{\varphi}\backslash\mathbb{R}$ if needed) insuring the uniqueness of the projection in Step 6.
\end{remark}

\begin{remark}\label{projHier}
 Concerning the sum in Step~6$_h$, a possible alternative, avoiding solving an algebraic system, would be to employ a hierarchical basis for $M_h^{\varphi}$ and $M_h^{\psi}$ and to use the $P^1$ component of $\varphi_h$ and $\psi_h$ instead of a $L^2$--projection. Such approach, if tempting considering computational costs, could result in important loss of information, negating the effort made in determining $\psi_h$ at Step~5.
\end{remark}

%Concerning the velocity correction, PROJECTION DU GRADIENT DE $\varphi$, composante tangente de la vitesse et reference Rannacher/Guermond.

\section{Numerical experiments}\label{simulation}

The goal of this section is to validate this family of methods numerically, illustrate the accuracy of some representatives and compare them.
%and finally demonstrate its numerical efficiency and advantages. 
Two tests will be presented: 
\begin{itemize}
    \item the first is an analytic test, based on a two dimensional manufactured solution presenting the time and spatial accuracy of four representive methods.
    \item The second test is the simulation of a flow past a cylinder of a Carreau fluid for various values of the power index~\cite{LasPraGia2012}. Here the \textit{drag} and \textit{wake} (or \textit{re-circulation}) lenght will be the object of comparison.
\end{itemize}
%a lid-driven cavity for an incompressible non--Newtonian fluid with a Carreau--Yasuda model for the viscosity. Although this test uses relatively high Reynolds number, it is found in numerous publications (see \cite{XXX,XXX} for some literature review on the lid-driven), and can be viewed as a "classical test". 
%
% This second test is found in numerous publications (see \cite{PatBhaChh2009, Pan2016} and the references therein, or \cite{MosJayMag2010} for a more complex law (Herschel-Bulkley)), and can be viewed as a "reference test".
%
% It uses more realistic data (boundary conditions, viscosity, etc.), has a well documented behaviour and offers reference values for punctual and boundary related physical quantities. %Moreover, as we imposes conditions such that a steady state is expected, it allows us to illustrates the behaviour of the method when used on phenomenon converging to a steady state. 
%
% As this test exhibit well documented behaviour for critical value of the power index, it is found in numerous publications, and can be viewed as a "reference test" (see \cite{PatBhaChh2009, Pan2016} and the references therein, or \cite{MosJayMag2010} for a more complex law (Herschel-Bulkley)).

\subsection{Choice of representatives and efficiency}
Section~\ref{Strategies} offers various approaches for the shear rate projection. These different strategies influences accuracy, efficiency and computational costs.
%For all the tests that follows we chose to retain both totally explicit and totally implicit schemes. 
It is not our intention to establish the relative merit of all those methods.
%
%In fact any of those methods can prove to be a valid choice, context of use must be considered for such comparison.
%
We will merely illustrate the capability of a small number of schemes when used on specific problems. For the first test we chose fully explicit and fully implicit representatives, leaving hybrid and Newton's variation to future, more thorough, analysis. For the second test, as we try to reproduce a steady state, the unsteady implicit shear rate projection is used and compared to a "steady" mixed algorithm and to results from the literature.

The use of the explicit strategy for both non linear terms makes \eqref{prediction3} and \eqref{AlgNSDivu_h} linear equations. This is tempting as it reduces computational efforts at each time step. However, totally explicit (as well as semi--explicit (hybrid) strategies) leads to conditional stability and certainly imposes limitations on the lenght of the time step. 
%For a simulation with a fixed final time, this restriction can prove to be costly from a numerical standpoint. 
Nevertheless it would be erroneous to exclude such approaches, as 
%in specific contexts (depending on the physical data) 
they can translate into very efficient algorithms.

Implicit strategies, fixed-point or Newton-like approaches, are of general use. 
%and can be qualified as less dependant on the physical data.
%
%Their 
Their efficiency is sometimes questioned since these choices leads to an iterative process. 
However implicit discretization minimizes the constraint on time step length making it possible to regain (when compared to explicit methods) numerical efficiency through the use of larger timestep (for a comparable precision). In those cases, it results in a better overall performance of the algorithm.

\begin{figure}[!htb]
\begin{center}
\includegraphics[width=0.475\textwidth]{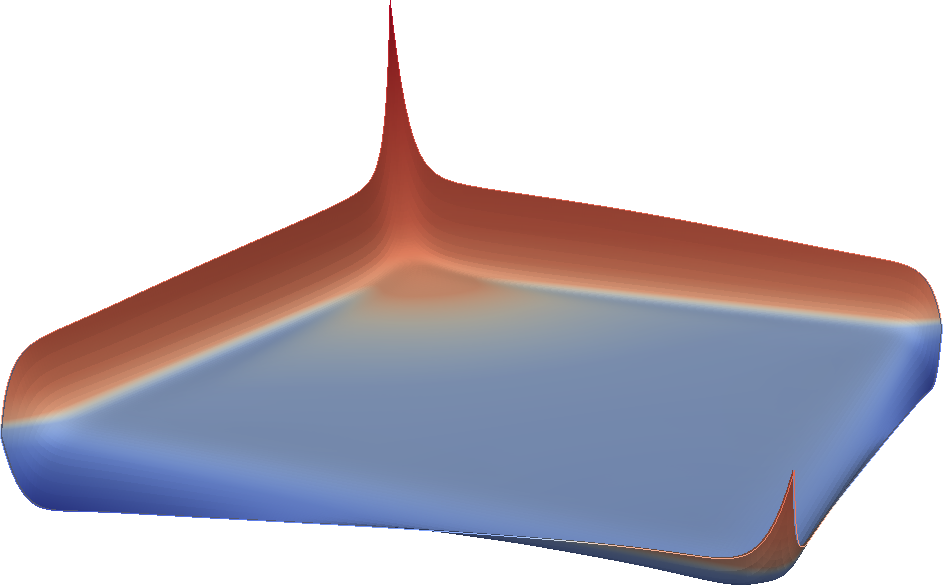}
\hspace*{0.03\textwidth}
\includegraphics[width=0.475\textwidth]{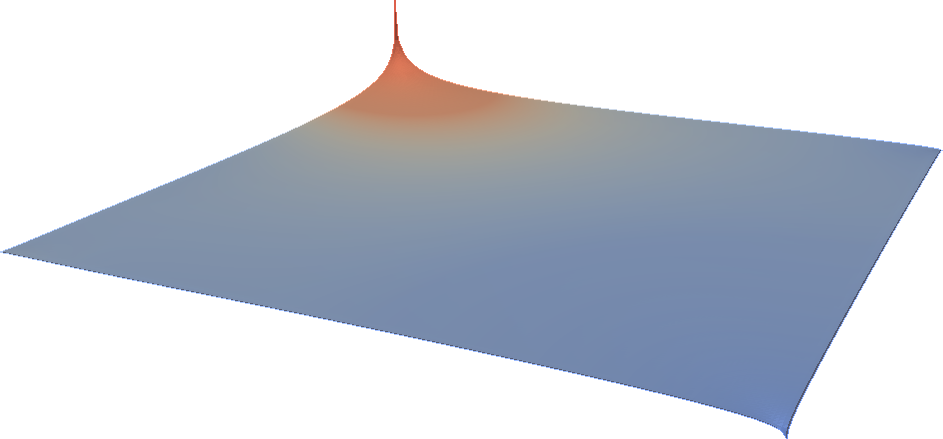}
% \qquad
% \includegraphics[width=0.42\textwidth]{fig2b.png}
\caption{Spatial distribution of the pressure error (signed, and rescale by a factor 10). On the left the implicit incremental projection method (IP$_{im}$), on the right the SRP$_{im}$ method. Data collected at time $t=1$ for $\Delta t = 1/80$ for a $200\times 200$ regular triangular mesh.}
\label{fig_para_pre}
\end{center}
\end{figure}

Four schemes will be used, 
\begin{enumerate}
\item a basic explicit incremental projection, denoted (IP$_{ex}$), obtained by imposing $\psi_h\equiv 0$,  ignoring Step~5$_h$ and putting,  at time $t^{n+1}$
$$\pVu^\star = \Vu^n,\qquad \nu^\star = \nu^{n+1}(\Vu^n),$$
\item a basic implicit incremental projection, denoted (IP$_{im}$), obtained by imposing $\psi_h\equiv 0$,  ignoring Step~5$_h$ and putting,  at time $t^{n+1}$
$$\pVu^\star = \pVu_k\qquad \nu^\star = \nu^{n+1}(\pVu_k),$$
\item the explicit shear rate projection scheme, denoted (SRP$_{ex}$), obtained by imposing, at time $t^{n+1}$
$$\pVu^\star = \Vu^n,\qquad \nu^\star = \nu^{n+1}(\Vu^n),$$
\item the implicit shear rate projection scheme, denoted (SRP$_{im}$), obtained by imposing, at iteration $k+1$ of time $t^{n+1}$
$$\pVu^\star = \pVu_k\qquad \nu^\star = \nu^{n+1}(\pVu_k).$$
\end{enumerate}
The algorithms are implemented in a finite element context using the FreeFem++ software, see \cite{freefem}. Note that all computations are performed using
%\begin{itemize}
%\item 
triangular mesh and Taylor-Hood interpolation: velocity are quadratic interpolation ($P_2$) and the pressure is linear ($P_1$). Insuring the respect of the inf--sup condition (see \cite{GueMinShe2006,BofBreFor2013}). Both projections, $\varphi$ and $\psi$, are taken in $P_2$.

\subsection{Accuracy test}
The first test is used to illustrate and compare the time and space accuracy of the shear rate projection (SRP) with the basic incremental projection (IP) for viscosity depending on the velocity of the fluid. As in \cite{DetYak2017} a variation on the finite element tests proposed by Guermond et al. in \cite{GueShe2004} is proposed

\begin{equation}\label{Sol_Analytique}
\left \{
\begin{array}{c}
u_1 =  \sin(x+t) \sin(y+t),\qquad %\vspace*{4pt}\\
u_2 = \cos(x+t)\cos(y+t)\vspace*{4pt}\\
p  =  \sin(x-y+t)\vspace*{4pt}\\
\rho \equiv 1,\,\qquad \nu_0 = 1\quad \nu(\Vs{u}) =  \nu_0(1+\|\Du\|^2)^{(m-1)/2}\quad m = 1/2\vspace*{4pt}\\
\Omega = ]0,1[\times]0,1[\\
\end{array}
\right.
\end{equation}
with a suitable forcing term depending on $\nu$.
% : 
% \begin{equation} \label{Test1_NS_F}
% \Vs{f} = \displaystyle \frac{ \partial \Vs{u}} {\partial t} \,+\, \left(\Vs{u} \cdot \nabla \right) \, \Vs{u}
% \, -\, \nabla \cdot \left(  2\nu(\Vs{u})\, \Du\right) \, + \, \nabla p 
% \end{equation}
\subsubsection{Time accuracy test}
We use a uniform mesh sufficiently fine, $200\times 200$, to insure negligible spatial error for the chosen range of time step. As for Step 2$_h$, when an implicit form is used, the stopping criterion of the fixed point uses the norm of the variations of $\pVu$ with a tolerance of $10^{-8}$. The velocity and pressure error to the exact solution will be computed using two types of norms, denoting $e$ the error, we will use
$$\|e\|_{\ell^2(S)} = \left(\Delta t \sum\limits_n \|e^n\|_S^2\right)^{1/2}\qquad \|e\|_{\ell^\infty(S)} = \max\limits_n \|e^n\|_S.$$

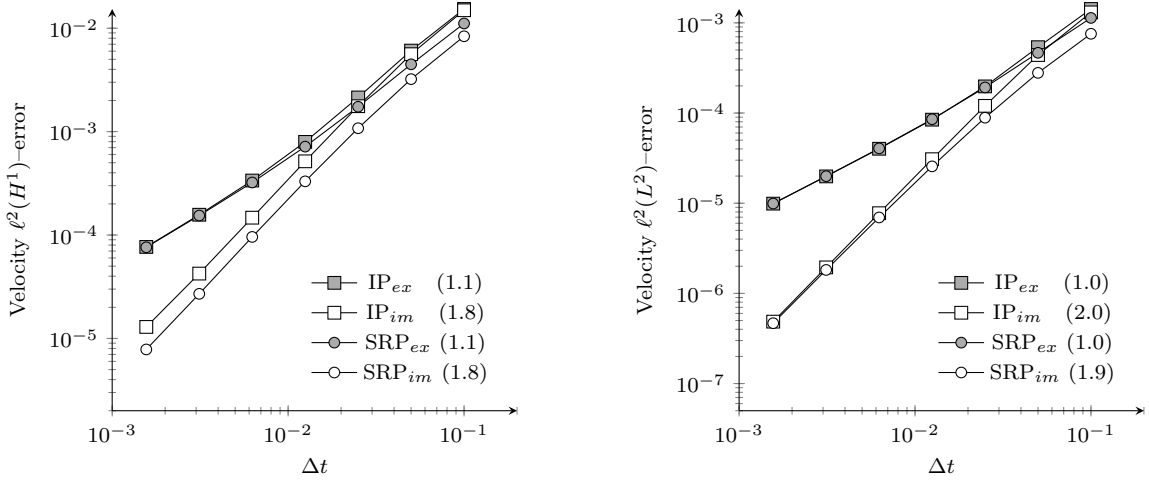
\begin{figure}[!ht]%[tpb]
\footnotesize
\begin{center}
%
% gauche
%
\begin{tikzpicture}[scale=1]
\begin{loglogaxis}[
width=.42\textwidth,
height=.42\textwidth,
xmax=0.2,
xmin=1e-3,
ymin=2e-6,
xlabel={$\Delta t$},
ylabel={Velocity $\ell^2(H^1$)--error},
axis lines=left,
legend style={draw=none, at={(.47,.37)},anchor=north west},
    cycle list name=black white]
\addplot [mark=square*, mark options={fill=black!33, scale=1.2}]
coordinates {
% (0.1       , 0.0153414  )
% (0.05      , 0.00607311 )
% (0.025     , 0.00213957 )
% (0.0125    , 0.000795833)
% (0.00625   , 0.000336438)
% (0.003125  , 0.000157099)
% (0.0015625 , 7.67416e-05)
(1./10  ,0.0153414)   
(1./20  ,0.00607311)  
(1./40  ,0.00213957)  
(1./80  ,0.000795833) 
(1./160  ,0.000336438)
(1./320  ,0.000157099)
(1./640  ,7.67416e-05)
};
\addlegendentry{IP$_{ex}\ \ $ $(1.1)$}

\addplot [mark=square*, mark options={fill=white, scale=1.2}]
coordinates {
(1./10  ,0.0149141)   
(1./20  ,0.00560859)  
(1./40  ,0.00176138)  
(1./80  ,0.0005162)   
(1./160  ,0.000147177)
(1./320  ,4.23849e-05)
(1./640  ,1.29222e-05)
};
\addlegendentry{IP$_{im}\ \ $ $(1.8)$}

\addplot [mark=*, mark options={fill=black!33}]
coordinates {
(1./10  ,0.0111098)   
(1./20  ,0.00447503)  
(1./40  ,0.00174071)  
(1./80  ,0.000715734) 
(1./160  ,0.000322116)
(1./320  ,0.000154386)
(1./640  ,7.61366e-05)
};
\addlegendentry{SRP$_{ex}$ $(1.1)$}

\addplot [mark=*, mark options={fill=white}]
coordinates {
(1./10  ,0.00835623)   
(1./20  ,0.00321756)   
(1./40  ,0.00107828)   
(1./80  ,0.000330423)  
(1./160  ,9.58718e-05) 
(1./320  ,2.69983e-05) 
(1./640  ,7.82744e-06)
};
\addlegendentry{SRP$_{im}$ $(1.8)$}
\end{loglogaxis}
\end{tikzpicture}
\hspace*{0.075\textwidth}
%
%                                                           droite
%
\begin{tikzpicture}[scale=1]
\begin{loglogaxis}[
width=.42\textwidth,
height=.42\textwidth,
xmax=0.2,
xmin=1.e-3,
ymin=5e-8,
xlabel={$\Delta t$},
ylabel={Velocity $\ell^2(L^2$)--error},
axis lines=left,
legend style={draw=none, at={(.47,.37)},anchor=north west},
    cycle list name=black white]
\addplot [mark=square*, mark options={fill=black!33, scale=1.2}]
coordinates {
% (0.1       , 0.00142465 )
% (0.05      , 0.000537056)
% (0.025     , 0.000197406)
% (0.0125    , 8.44013e-05)
% (0.00625   , 4.01634e-05)
% (0.003125  , 1.98351e-05)
% (0.0015625 , 9.89152e-06)
(1./10  ,0.00142465)  
(1./20  ,0.000537056) 
(1./40  ,0.000197406) 
(1./80  ,8.44013e-05) 
(1./160  ,4.01634e-05)
(1./320  ,1.98351e-05)
(1./640  ,9.89152e-06) 
};
\addlegendentry{IP$_{ex}\ \ $ $(1.0)$}
\addplot [mark=square*, mark options={fill=white, scale=1.2}]
coordinates {
(1./10  ,0.00130613)  
(1./20  ,0.000439785) 
(1./40  ,0.00012006)  
(1./80  ,3.07628e-05) 
(1./160  ,7.74938e-06)
(1./320  ,1.94259e-06)
(1./640  ,4.86213e-07)
};
\addlegendentry{IP$_{im}\ \ $ $(2.0)$}
% %
\addplot [mark=*, mark options={fill=black!33}]
coordinates {
(1./10  ,0.00113376)  
(1./20  ,0.000465466) 
(1./40  ,0.00019218)  
(1./80  ,8.53415e-05) 
(1./160  ,4.05659e-05)
(1./320  ,1.99371e-05)
(1./640  ,9.91546e-06)
};
\addlegendentry{SRP$_{ex}$ $(1.0)$}
% %
\addplot [mark=*, mark options={fill=white}]
coordinates {
(1./10  ,0.000756258)  
(1./20  ,0.0002789)    
(1./40  ,8.87916e-05)  
(1./80  ,2.56123e-05)  
(1./160  ,6.95078e-06) 
(1./320  ,1.82033e-06) 
(1./640  ,4.67281e-07)
};
\addlegendentry{SRP$_{im}$ $(1.9)$}
\end{loglogaxis}
\end{tikzpicture}
\caption{ Velocity error for the basic incremental projection explicit and implicit form (IP$_{ex}$ and IP$_{im}$) and the shear rate projection explicit and implicit method (SRP$_{ex}$ and SRP$_{im}$). On the left, the $H^1$ norm  and on the right the $L^2$ norm of the error. In parenthesis the approximated slope, data collected at time $t=1$ for $\Delta t = 1/10$ to $ 1/640$ on a $200\times 200$ regular triangular mesh.}
\label{fig_acc_time_velo}
\end{center}
\end{figure}

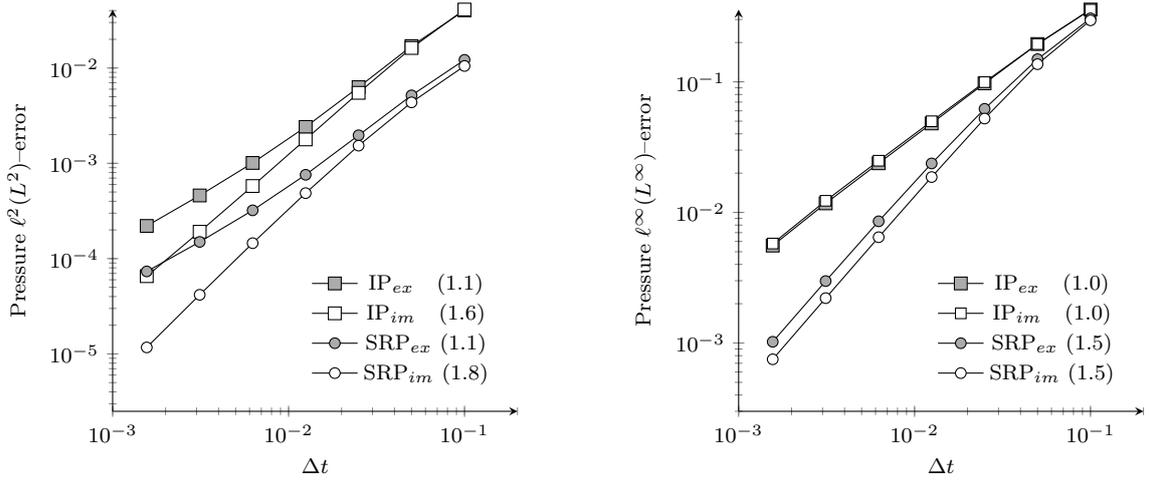
\begin{figure}[!ht]
\footnotesize
\begin{center}
\begin{tikzpicture}[scale=1]
\begin{loglogaxis}[
width=.42\textwidth,
height=.42\textwidth,
xmax=0.2,
xmin=1.e-3,
ymin=2.5e-6,
xlabel={$\Delta t$},
ylabel={Pressure $\ell^2(L^2)$--error},
axis lines=left,
legend style={draw=none, at={(.47,.37)},anchor=north west},
    cycle list name=black white]
\addplot [mark=square*, mark options={fill=black!33, scale=1.2}]
coordinates {
(1./10  ,0.0402177)   
(1./20  ,0.016988)   (
(1./40  ,0.00629447)  
(1./80  ,0.00240772)  
(1./160  ,0.00100925) 
(1./320  ,0.000459995)
(1./640  ,0.000220248) 
};
\addlegendentry{IP$_{ex}\ \ $ $(1.1)$}
\addplot [mark=square*, mark options={fill=white, scale=1.2}]
coordinates {
(1./10  ,0.0410372)   
(1./20  ,0.0162291)   
(1./40  ,0.00548929)  
(1./80  ,0.00178183)  
(1./160  ,0.000577673)
(1./320  ,0.000191408)
(1./640  ,6.54859e-05)
};
\addlegendentry{IP$_{im}\ \ $ $(1.6)$}
% %
\addplot [mark=*, mark options={fill=black!33}]
coordinates {
(1./10  ,0.0121738)   
(1./20  ,0.00516455)  
(1./40  ,0.00196925)  
(1./80  ,0.00075814)  
(1./160  ,0.000321086)
(1./320  ,0.00014975) 
(1./640  ,7.34187e-05)
};
\addlegendentry{SRP$_{ex}$ $(1.1)$}
% %
\addplot [mark=*, mark options={fill=white}]
coordinates {
(1./10  ,0.0105037)   
(1./20  ,0.00436133)  
(1./40  ,0.0015327)   
(1./80  ,0.00048652)  
(1./160  ,0.000145221)
(1./320  ,4.17205e-05)
(1./640  ,1.17031e-05)
};
\addlegendentry{SRP$_{im}$ $(1.8)$}
\end{loglogaxis}
\end{tikzpicture}
\hspace*{0.075\textwidth}
%
%                                                          L^\infty
%
\begin{tikzpicture}[scale=1]
\begin{loglogaxis}[
width=.42\textwidth,
height=.42\textwidth,
xmax=0.2,
xmin=1.e-3,
ymin=3e-4,
xlabel={$\Delta t$},
ylabel={Pressure $\ell^\infty(L^\infty)$--error},
axis lines=left,
legend style={draw=none, at={(.47,.37)},anchor=north west},
    cycle list name=black white]
\addplot [mark=square*, mark options={fill=black!33, scale=1.2}]
coordinates {
(1./10  ,0.357664)
(1./20  ,0.194831)
(1./40  ,0.0973274)
(1./80  ,0.0480235)
(1./160  ,0.0237976)
(1./320  ,0.0117333)
(1./640  ,0.00557253)
};
\addlegendentry{IP$_{ex}\ \ $ $(1.0)$}
\addplot [mark=square*, mark options={fill=white, scale=1.}]
coordinates {
(1./10  ,0.357241)
(1./20  ,0.194496)
(1./40  ,0.099454)
(1./80  ,0.0499102)
(1./160  ,0.0248825)
(1./320  ,0.0122804)
(1./640  ,0.0057786)
};
\addlegendentry{IP$_{im}\ \ $ $(1.0)$}
% %
\addplot [mark=*, mark options={fill=black!33}]
coordinates {
(1./10  ,0.306262)
(1./20  ,0.149095)
(1./40  ,0.0621346)
(1./80  ,0.0237264)
(1./160  ,0.00855269)
(1./320  ,0.00298153)
(1./640  ,0.00102182)
};
\addlegendentry{SRP$_{ex}$ $(1.5)$}
% %
\addplot [mark=*, mark options={fill=white}]
coordinates {
(1./10  ,0.29579)
(1./20  ,0.136514)
(1./40  ,0.0523397)
(1./80  ,0.0186623)
(1./160  ,0.00646)
(1./320  ,0.0022072)
(1./640  ,0.000750611)
};
\addlegendentry{SRP$_{im}$ $(1.5)$}
\end{loglogaxis}
\end{tikzpicture}
\caption{ Pressure error ($L^2$ and $L^\infty$ norm) for the basic incremental projection explicit and implicit form (IP$_{ex}$ and IP$_{im}$) and the shear rate projection explicit and implicit method (SRP$_{ex}$ and SRP$_{im}$). In parenthesis the approximated slope, data collected at time $t=1$ for $\Delta t = 1/10$ to $ 1/640$ on a $200\times 200$ regular triangular mesh.}
\label{fig_acc_time_pre}
\end{center}
\end{figure}

As expected, see  Figure~\ref{fig_acc_time_velo}--\ref{fig_acc_time_pre}, the explicit form of both methods, IP$_{ex}$ and SRP$_{ex}$, exhibit a loss of accuracy when compared to their implicit version. The pressure error, Figure~\ref{fig_para_pre} and \ref{fig_acc_time_pre} illustrate the gain in accuracy of this "extended" SRP method, this is comparable to the results presented in \cite{DetYak2017}.

% \newpage
%
%
%                                                            Figure Vitesse
%
%
%
%                                              Figure Pression
%
% \clearpage
For the velocity in $L^2$--norm, (Figure~\ref{fig_acc_time_velo} right), both explicit methods exhibit a linear rate while both implicit methods have a quadratic convergence rate. The SRP does not produce any significant gain in precision for the velocity. The $H^1$ convergence rates are identical: the explicit methods have a rate of $1.1$ while the implicit method shows a rate of almost $2$, this is explained by the increased accuracy of the pressure. The SRP$_{im}$ produces a velocity  error in $H^1$ norm approximately $2$ time smaller then the $H^1$ error produced by the implicit incremental projection.

%In the present case, the 
The $L^2$ rates for the pressure (Figure~\ref{fig_acc_time_pre} left) seems relatively similar for both explicit and implicit pair of schemes (around 1 for the explicit schemes and 1.7 for the implicit ones). However both version of the SRP are clearly more accurate, reducing the $L^2$--error by a factor roughly equal to $3$. From Figure~\ref{fig_acc_time_pre} the $L^\infty$ rate of both IP methods is linear and the SRP's rate are $3/2$: the extreme values of the errors of the SRP are decaying more rapidly then those of the IP method. Figure~\ref{fig_para_pre} illustrate at time $t=1$, for $\Delta t = 1/80$, the spatial distribution of the error for the pressure; giving a sense of the impact of $\psi$ on the pressure.

In conclusion, for non-homogeneous viscosity, both version of the SRP method produce a more accurate approximation of the pressure (and velocity) than the IP version, even in the explicit case, the use of the SRP improves the pressure approximation. 

\subsubsection{Spatial accuracy test}
Since the SRP schemes contain an additional problem to solve, \eqref{EqPsi_h}, it is necessary to illustrate the impact of this added computations on the spatial accuracy of the scheme. Here we included two other methods as comparative: IP$_{ex}$ which is the simplest transition to projection methods for generalized Newtonian fluids and which has the smallest computational cost per time step of all the methods presented and an implicit version of the very simple penalty method~\cite{Tem1977}, also called mixed method, using a penalty of $10^{-10}$ for the incompressibility.
%%%%%%%%%%%%%%%%%%%%
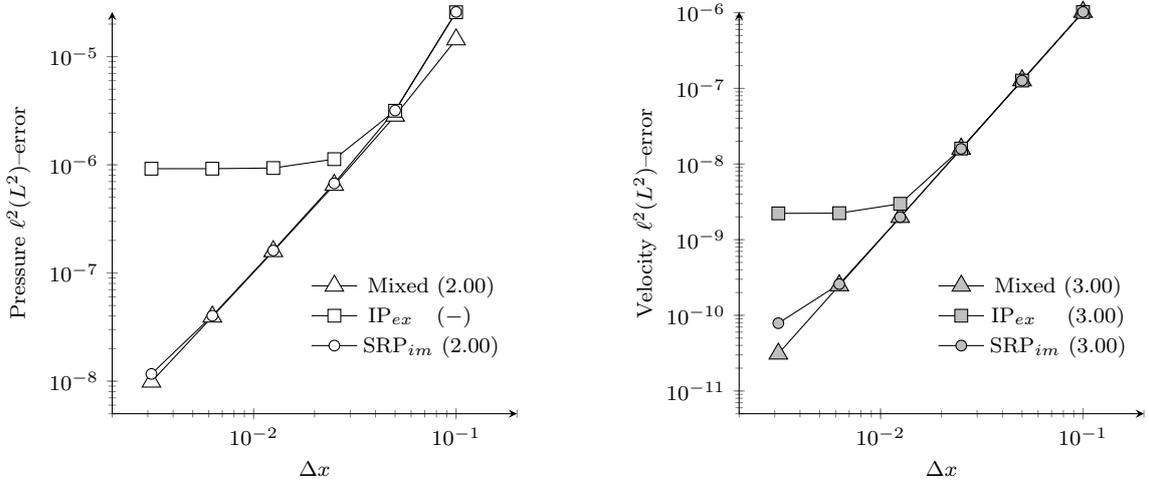
\begin{figure}[hbt]
\footnotesize
\begin{center}
\begin{tikzpicture}[scale=1]
\begin{loglogaxis}[
width=.42\textwidth,
height=.42\textwidth,
xmax=1/5,
xmin=1/500,
ymin=5e-9,
xlabel={$\Delta x$},
ylabel={Pressure $\ell^2(L^2)$--error},
axis lines=left,
legend style={draw=none, at={(.47,.37)},anchor=north west},
    cycle list name=black white]

\addplot [mark=triangle*, mark options={fill=white, scale = 2}]
coordinates {
(1./10  ,1.4396e-05)   
(1./20  ,2.83784e-06)  
(1./40  ,6.50714e-07)  
(1./80  ,1.58691e-07)  
 (1./160  ,3.94212e-08)
 (1./320  ,9.86049e-09)
};
% \addlegendentry{$\text{Mixed : }\|p-p_h\|_{\ell^2(L^2)}$ $(2.00)$}
\addlegendentry{$\text{Mixed }(2.00)$}

\addplot  [mark=square*, mark options={fill=white, scale = 1.2}]%[mark = round,mark options={scale=.75}] 
coordinates {
(1./10  ,2.57065e-05) 
(1./20  ,3.16101e-06) 
(1./40  ,1.12816e-06) 
(1./80  ,9.33794e-07) 
(1./160  ,9.21056e-07)
(1./320  ,9.20247e-07)
 };
% \addlegendentry{$\text{IP$_{ex}$\quad : }\|p-p_h\|_{\ell^2(L^2)}$ $(-)\ \ \ $}
\addlegendentry{$\text{IP$_{ex}$  }\ \ (-)\ \ \ $}

 \addplot  [mark=*, mark options={fill=white}]
 coordinates {
% (1/ 20 , 8.76336e-05  )
% (1/ 40 , 9.93268e-06  )
% (1/ 80 , 2.27944e-06  )
% (1/160 , 5.61741e-07  )
% (1/320 , 1.44778e-07 )
(1./10  ,2.58811e-05)   
(1./20  ,3.17947e-06)   
(1./40  ,6.73694e-07)   
(1./80  ,1.61486e-07)   
  (1./160  ,4.02475e-08)
  (1./320  ,1.16538e-08)
 };
%  \addlegendentry{$\text{SRP$_{im}$: }\|p-p_h\|_{\ell^2(L^2)}$ $(2.00)$}
 \addlegendentry{$\text{SRP$_{im}$ }(2.00)$}

% \addlegendimage{empty legend}
% \addlegendentry{}
% \addlegendimage{empty legend}
% \addlegendentry{}
\end{loglogaxis}
\end{tikzpicture}
\hspace*{0.075\textwidth}
%
%                                                          velo
%
\begin{tikzpicture}[scale=1]
\begin{loglogaxis}[
width=.42\textwidth,
height=.42\textwidth,
xmax=1/5,
xmin=1/500,
ymin=5e-12,
xlabel={$\Delta x$},
ylabel={Velocity $\ell^2(L^2$)--error},
axis lines=left,
legend style={draw=none, at={(.47,.37)},anchor=north west},
    cycle list name=black white]
%
% \addplot [mark=triangle*, mark options={fill=black!25, scale=2}]
% coordinates {
% (1./10  ,7.72491e-05) 
% (1./20  ,1.93259e-05) 
% (1./40  ,4.83131e-06) 
% (1./80  ,1.20779e-06) 
% (1./160  ,3.01944e-07)
% (1./320  ,7.54857e-08)
% };
% \addlegendentry{$\text{Mixed  : }\|u-u_h\|_{\ell^2(H^1)}$ $(2.00)$}

% \addplot [mark=square*, mark options={fill=black!25, scale = 1.2}]
% coordinates {
% (1./10  ,7.70785e-05) 
% (1./20  ,1.92561e-05) 
% (1./40  ,4.82121e-06) 
% (1./80  ,1.20676e-06) 
% (1./160  ,3.02982e-07)
% (1./320  ,8.01367e-08)
% };
% \addlegendentry{$\text{IP$_{ex}$\quad   : }\|u-u_h\|_{\ell^2(H^1)}$ $(2.00)$}

% \addplot [mark=*, mark options={fill=black!25}]
% coordinates {
% % (1/ 20 , 0.00027296   )
% % (1/ 40 , 6.80569e-05 )
% % (1/ 80 , 1.70278e-05 )
% % (1/160 , 4.2592e-06 )
% % (1/320 , 1.06552e-06 )
% (1./10  ,7.70596e-05)  
% (1./20  ,1.92546e-05)  
% (1./40  ,4.821e-06)    
% (1./80  ,1.20644e-06)  
%  (1./160  ,3.01786e-07)
%  (1./320  ,7.55187e-08)
% };
% \addlegendentry{$\text{SRP$_{im}$: }\|u-u_h\|_{\ell^2(H^1)}$ $(2.00)$}

% \addlegendimage{empty legend}
% \addlegendentry{}

\addplot [mark=triangle*, mark options={fill=black!25, scale=2.}]
coordinates {
(1./10  ,1.02342e-06)
(1./20  ,1.27404e-07)
(1./40  ,1.5899e-08) 
(1./80  ,1.9864e-09) 
(1./160  ,2.4827e-10)
(1./320  ,3.1047e-11)
};
\addlegendentry{$\text{Mixed  }(3.00)$}

\addplot [mark=square*, mark options={fill=black!25, scale = 1.2}]
coordinates {
(1./10  ,1.023e-06)   
(1./20  ,1.26934e-07) 
(1./40  ,1.60162e-08) 
(1./80  ,2.98081e-09) 
(1./160  ,2.23933e-09)
(1./320  ,2.22583e-09)
};
\addlegendentry{$\text{IP$_{ex}$ }\ \ \ (3.00)$}

\addplot [mark=*, mark options={fill=black!25}]
coordinates {
% (1/ 20 ,1.80092e-06)
% (1/ 40 , 2.23959e-07 )
% (1/ 80 , 2.80026e-08 )
% (1/160 , 3.50192e-09 )
%  (1/320 , 4.3815e-10 )
(1./10  ,1.02259e-06)  
(1./20  ,1.26878e-07)  
(1./40  ,1.58601e-08)  
(1./80  ,1.98524e-09)  
(1./160  ,2.58408e-10) 
(1./320  ,7.85926e-11) 
};
\addlegendentry{$\text{SRP$_{im}$ }(3.00)$}
\end{loglogaxis}
\end{tikzpicture}
\caption{Spatial accuracy for the mixed method (penalized incompressibility), explicit incremental method (IP$_{ex}$) and implicit shear rate projection (SRP$_{im}$). Pressure (left) and velocity (right) error to the analytical solution after 50 time steps of lenght $\Delta t = 10^{-4}$. Data collected for regular triangular meshes of size $\Delta x = 1/10$ to $1/320$.}
\label{fig_ordre_espace}
\end{center}
\end{figure}

For this accuracy test, 6 meshes obtained by regular subdivision from $1/10$ to $1/320$ were used. A total of 50 time steps of fixed lenght $10^{-4}$ (for a final time t=0.005) was used; assuming such small time step would limit the impact of the time discretization on the approximation error. 
The first observation concerns the IP$_{ex}$ scheme, the $L^2$--error for the velocity and pressure seems to be rapidly saturated, mesh refinement over $1/80$ having no effect on the error for the current time step of $10^{-4}$. 
%\JDm{Que dire de la norme H1 ? la norme L2 est petite, H1 = semi-H1 a peu près, et alors?.} 
From the finest mesh, the SRP$_{im}$ seems to exhibit a faster saturation of the $L^2$--error compared to the penalty method (denoted Mixed in Figure~\ref{fig_ordre_espace}).
Nevertheless, the SRP$_{im}$ and Mixed method produce solutions with comparable spatial convergence rate, confirming the good behaviour in space and in time of $\psi_h$ solution of \eqref{EqPsi_h}.
% \newpage
\subsection{Flow past a cylinder}
\begin{figure}[ht]
\centering
\begin{tikzpicture}
\tikzstyle{fleche}=[->,>=latex]
\tikzset{
           xmin/.store in=\xmin, xmin/.default=0., xmin=0.,
           xmax/.store in=\xmax, xmax/.default=4.75, xmax=5.1,
           ymin/.store in=\ymin, ymin/.default=0., ymin=0.,
           ymax/.store in=\ymax, ymax/.default=4.75, ymax=5.1,
}

 \newcommand {\fenetre}
 {\clip (\xmin-1.5,\ymin) rectangle (\xmax+1.5,\ymax);}
 \fenetre

\draw[] (\xmin,\ymin) rectangle (\xmax-0.35,\ymax-0.35); 
% \draw[](3.5, 4.5)      node {$\Omega$};
\draw[dotted, very thin] (2.38,\ymax-0.35) -- (2.38,\ymin);
\draw[dotted, very thin] (\xmin,2.38) -- (\xmax-0.35,2.38);

\draw[] (\xmin+0.3,2.38)  node  {$\Gamma_{i}$}; 
\draw[] (\xmax-0.65,2.38)  node  {$\Gamma_{o}$}; 
\draw[] (2.38,\ymax-0.65)  node  {$\Gamma_{t}$}; 
\draw[] (2.38,\ymin+0.3)  node  {$\Gamma_{t}$}; 
\draw[] (2.98,1.4)  node  {$\Gamma_{c}$}; 
\draw[] (1,1)  node  {$\Omega$}; 

\draw[](1.3, \ymax-0.17)  node  {$\mathnormal{\ell_{u}}$};
\draw[] (\xmin,\ymax-0.35) -- (\xmin,\ymax);
\draw[] (2.38,\ymax-0.35) -- (2.38,\ymax);
\draw[] (\xmax-0.35,\ymax-0.35) -- (\xmax-0.35,\ymax);
\draw[fleche, dashed, very thin] (1.5,\ymax-0.17) -- (2.38,\ymax-0.17);
\draw[fleche, dashed, very thin] (1.1,\ymax-0.17) -- (\xmin,\ymax-0.17);

\draw[](3.3, \ymax-0.17)  node  {$\mathnormal{\ell_{d}}$};
\draw[] (\xmax-0.35,\ymin) -- (\xmax,\ymin);
\draw[] (\xmax-0.35,\ymax-0.35) -- (\xmax,\ymax-0.35);
\draw[fleche, dashed, very thin] (3.1,\ymax-0.17) -- (2.38,\ymax-0.17);
\draw[fleche, dashed, very thin] (3.5,\ymax-0.17) -- (\xmax-0.35,\ymax-0.17);

\draw[](\xmax-0.15, 2.38)    node {$\mathnormal{h}$};
\draw[fleche, dashed, very thin] (\xmax-0.17,2.58) -- (\xmax-0.17,\ymax-0.35);
\draw[fleche, dashed, very thin] (\xmax-0.17,2.18) -- (\xmax-0.17,\ymin);

\draw[black] [fill=white](2.38,2.38) circle (0.2cm);
\draw[](2.38, 3.1)  node {$\mathnormal{D}$};
\draw[very thin] (2.175,2.475) -- (2.175,2.98);
\draw[very thin] (2.585,2.475) -- (2.585,2.98);
\draw[fleche, dashed, very thin] (1.88,2.85) -- (2.18,2.85);
\draw[fleche, dashed, very thin] (2.88,2.85) -- (2.58,2.85);

\draw[] (\xmin-0.9,2.38)  node  {$U_{0}$}; 
\draw[fleche, very thin] (\xmin-0.5,\ymin+0.1) -- (\xmin,\ymin+0.1);
\draw[fleche, very thin] (\xmin-0.5,\ymin+0.4) -- (\xmin,\ymin+0.4);
\draw[fleche, very thin] (\xmin-0.5,\ymin+0.7) -- (\xmin,\ymin+0.7);
\draw[fleche, very thin] (\xmin-0.5,\ymin+1.) -- (\xmin,\ymin+1.);
\draw[fleche, very thin] (\xmin-0.5,\ymin+1.3) -- (\xmin,\ymin+1.3);
\draw[fleche, very thin] (\xmin-0.5,\ymin+1.6) -- (\xmin,\ymin+1.6);
\draw[fleche, very thin] (\xmin-0.5,\ymin+1.9) -- (\xmin,\ymin+1.9);
\draw[fleche, very thin] (\xmin-0.5,\ymin+2.2) -- (\xmin,\ymin+2.2);
\draw[fleche, very thin] (\xmin-0.5,\ymin+2.5) -- (\xmin,\ymin+2.5);
\draw[fleche, very thin] (\xmin-0.5,\ymin+2.8) -- (\xmin,\ymin+2.8);
\draw[fleche, very thin] (\xmin-0.5,\ymin+3.1) -- (\xmin,\ymin+3.1);
\draw[fleche, very thin] (\xmin-0.5,\ymin+3.4) -- (\xmin,\ymin+3.4);
\draw[fleche, very thin] (\xmin-0.5,\ymin+3.7) -- (\xmin,\ymin+3.7);
\draw[fleche, very thin] (\xmin-0.5,\ymin+4.0) -- (\xmin,\ymin+4.0);
\draw[fleche, very thin] (\xmin-0.5,\ymin+4.3) -- (\xmin,\ymin+4.3);
\draw[fleche, very thin] (\xmin-0.5,\ymin+4.65) -- (\xmin,\ymin+4.65);

\end{tikzpicture}
\caption{Schematic of the domain. The axe of the cylinder correspond to the origin $(0,0)$.}
\label{fig_cyl}
\end{figure}
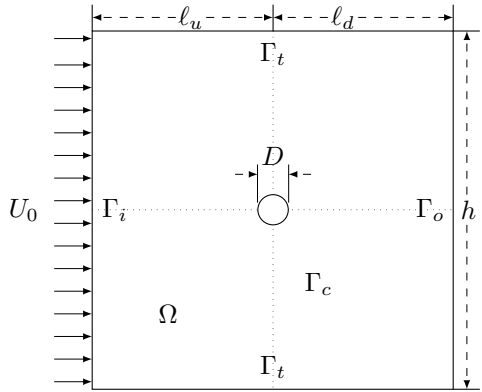

For this test, a two-dimensional model of the flow of a generalized Newtonian (Carreau) fluid past a cylinder, see figure~\ref{fig_cyl}, is considered.
This second test uses more realistic data (boundary conditions, viscosity, etc.) and can be viewed as a "reference test" as it is presented frequently (see \cite{LasPraGia2012,PatBhaChh2009, Pan2016} and the references therein, or \cite{MosJayMag2010} for a more complex law (Herschel-Bulkley)).

%%%%%%%%%%%%%%%%%%%%%%%%%%%%

This well documented test offers reference values for punctual and boundary related physical measures. This allows us to illustrate the accuracy of the projections schemes based on relevant "local" physical quantities of interest: the drag coefficient $C_D$ and the wake lenght $L$,
$$C_D = -\frac{2}{\rho}\int_{\Gamma_c} \sigma(\Vs{u})\Vs{n}\cdot \Vs{n}\ d\Gamma,\qquad L = \argmin\limits_{x> D/2}|\Vs{u}(x,0)|.$$
The fluid is characterized by the Carreau law, defining the viscosity, and by two adimensional numbers, the Renolds number ($Re$) and the Carreau number $C_U$, 
$$\nu(\Vs{u}) =  \nu_\infty + (\nu_0-\nu_\infty)\left(1+2(\lambda\|\Du\|)^2\right)^{(m-1)/2},\quad Re= \frac{\rho D U_0}{\nu_0},\quad C_U = \frac{\lambda U_0}{D}.$$
For this test we are interested in a steady state flow, which correspond to a relatively small $Re$. Contrary to tests frequently encountered in the literature, we will limit our tests to a single Reynolds number.
Experimentation will be made using the implicit shear rate projection (showing the ability of the SRP$_{im}$ to reproduce steady states results). 
Reference solutions will be computed using a steady state mixed algorithm (denoted SS-Mix) based on a penalty method inside a fixed point loop for the viscosity. 

The data for the tests are taken from~\cite{LasPraGia2012,Pan2016}, adding another point of comparison. 
%
%\subsubsection{Details of the numerical experiment}
%
%Although its physical relevance can be questioned, following~\cite{Pan2016} we treat the case of unconfined flow past a cylinder at Reynolds number $Re=10$. 
%It is well known that the boundary of the computational domain can influence the results if too close to the obstacle, this influence increase as the Reynold number is small. This explains the dimension The dimension of the domain will insure almost inexistant boundary effect on the results presented here.
We can summarize the parameters as follow:
\begin{equation}\label{donnees_ecoul_cyl}
\left \{
\begin{array}{c}
% \nu(\Vs{u}) =  \nu_\infty + (\nu_0-\nu_\infty)\left(1+(\lambda\|\Du\|)^2\right)^{(m-1)/2}\\
\Vu = (U_0,0)\quad\text{on }\Gamma_i,\quad U_0 = 1,\\
\Vu = (0,0)\quad\text{on }\Gamma_d,\quad \Vu\cdot \Vn = 0\quad\text{on }\Gamma_t.\\
\nu_0 = 1,\quad\nu_\infty = 0.001,\\
m\in\{0.4, 0.5, 0.6, 0.7, 0.8, 0.9, 1.0\},\\
\lambda = C_{U}\in \{10,20 \},\quad\rho = Re =10,\\
D= 1,\quad \ell_{u} = \ell_d = h/2,\quad h = 50000D.
\end{array}
\right.
\end{equation}
Finally, contrary to ~\cite{LasPraGia2012, PatBhaChh2009}, a natural (or open) boundary condition is applied on $\Gamma_o$. 

As an unconfined flow is modeled, and considering the open boundary condition, the outer boundary of the domain need to be sufficiently far from the cylinder so that boundary effects are at a minimum. This implies a relatively large domain. 
It should be notice that the dimensions chosen guarantees us of eliminating any boundary effects; however the use of a smaller domain could be possible (provided a proper domain dependency analysis is made). 
Mesh adaptation is used to reduce the number of elements by refining and de-refining locally the mesh, this makes it possible to have reasonable computational time even for such large domain.
For each power index value ($m$), a specifically adapted mesh will be used. These meshes are constructed based on the SS-Mix solutions for each value of $m$ (see~\cite{freefem} for details concerning the mesh adaptation procedure).

From~\cite{LasPraGia2012}, the critical Reynolds, at which the flow becomes unsteady, for this Carreau fluid and for a power index $m\ge 0.4$ is greater then $10$. Therefore the flows simulated here are certainly steady for all the values of $m$ used, the steady state will be reached using the SRP$_{im}$ scheme, an unsteady scheme. The simulation is stopped once the steady state is attained, as indicated by the stationnarity of certain physical quantities (here the drag coefficient and magnitude of the velocity).

\begin{figure}[!htb]
\centering
\begin{tabular}{cc}
\includegraphics[width=0.425\textwidth]{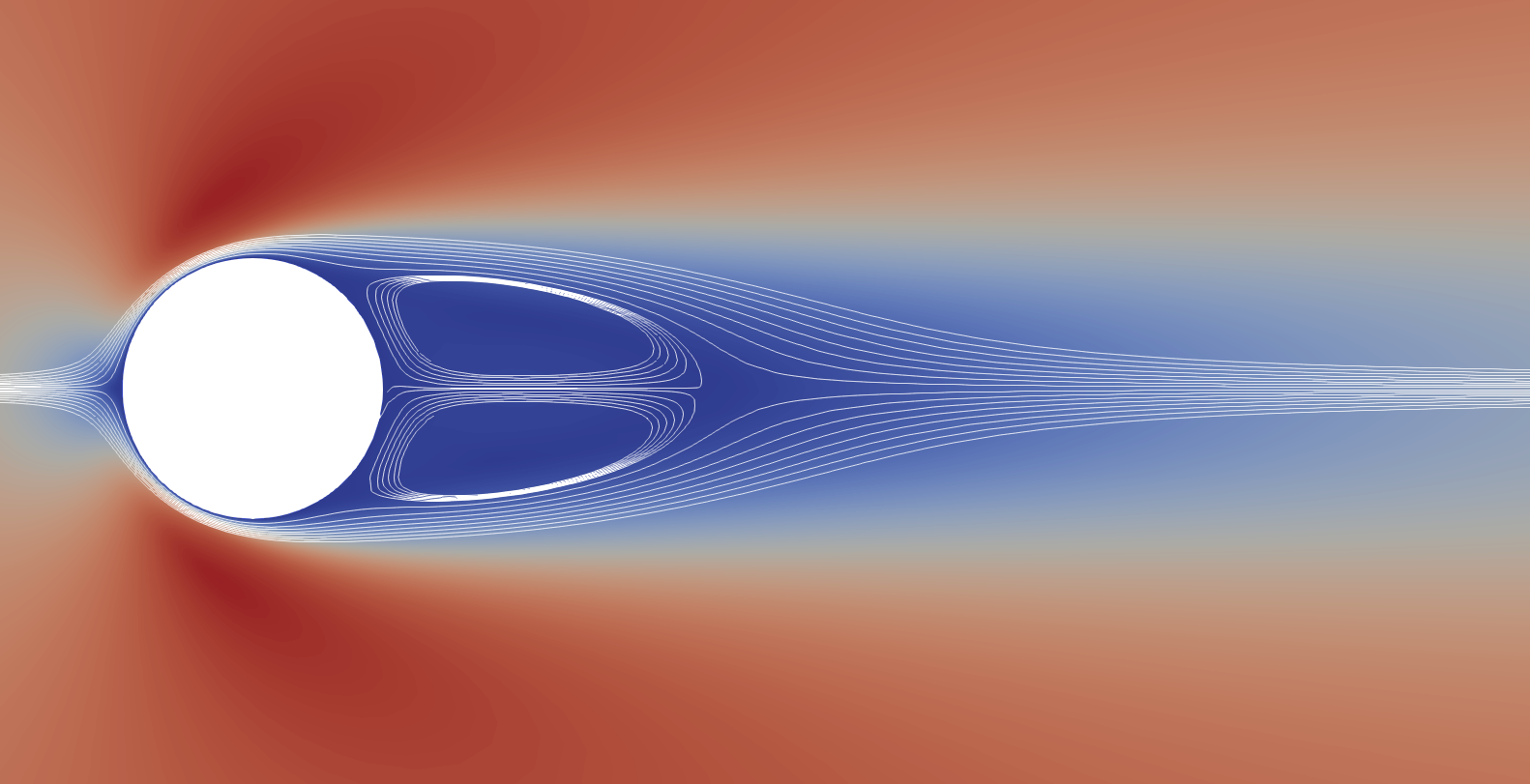} &\quad
\includegraphics[width=0.425\textwidth]{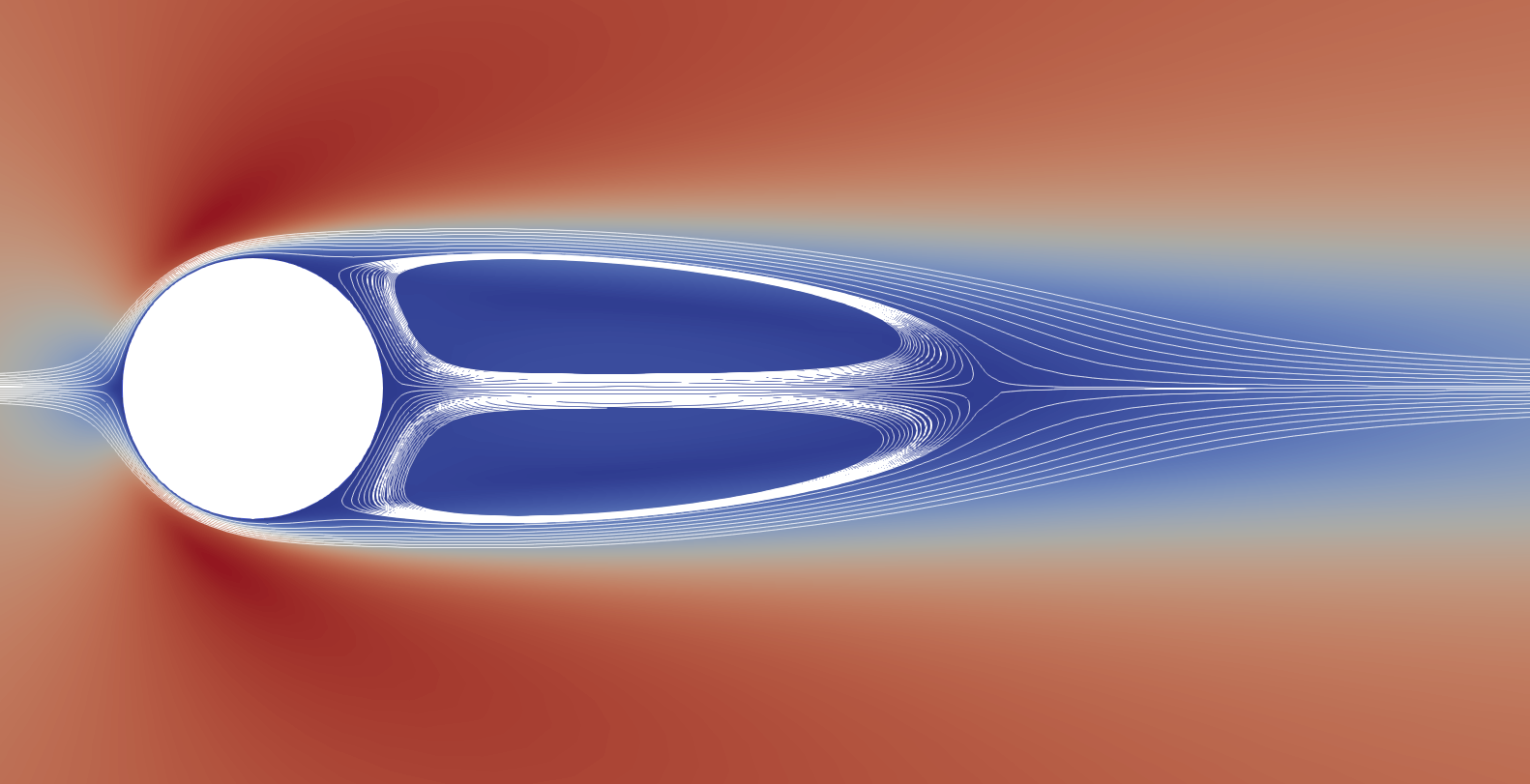}\\
\includegraphics[width=0.425\textwidth]{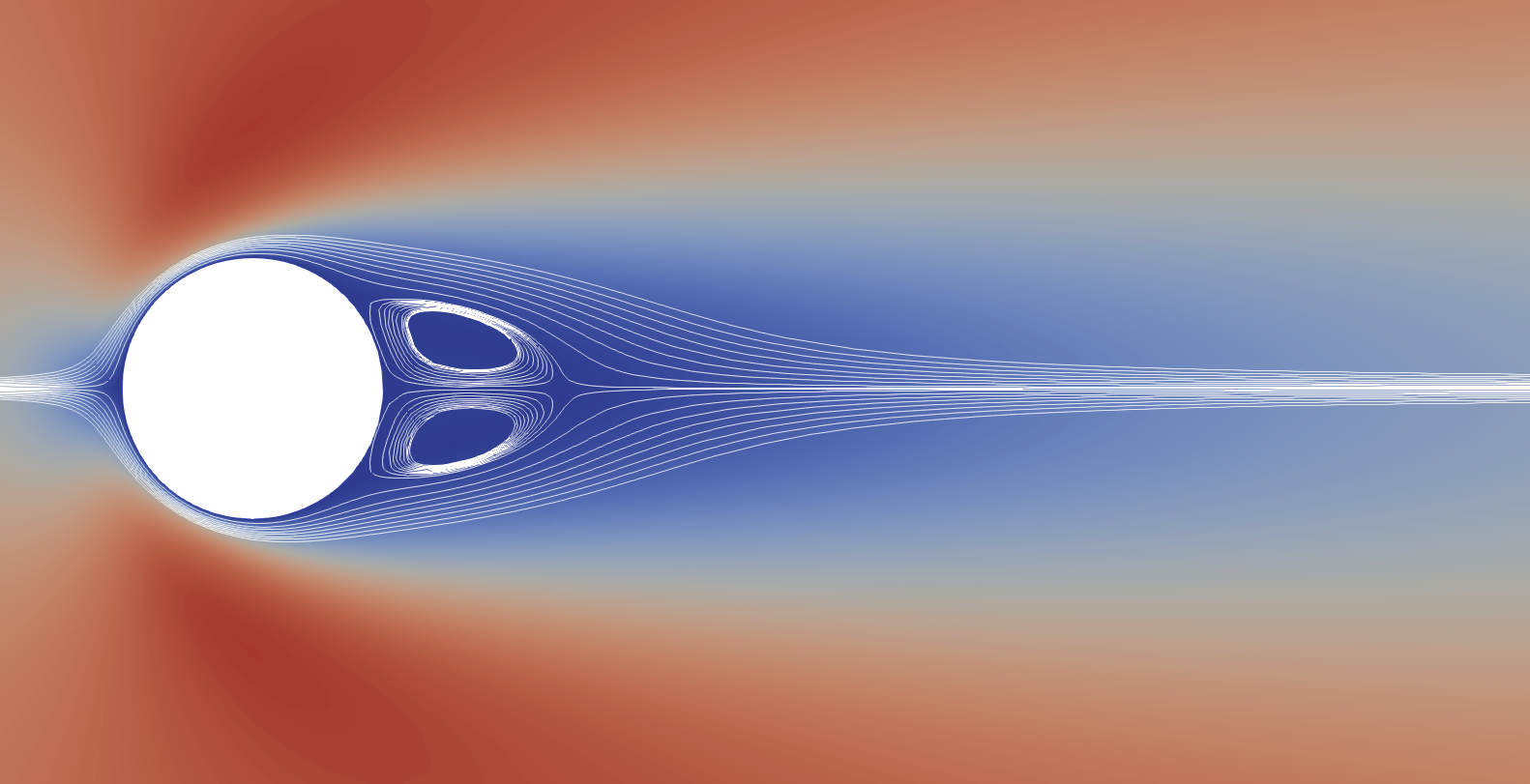} &\quad
\includegraphics[width=0.425\textwidth]{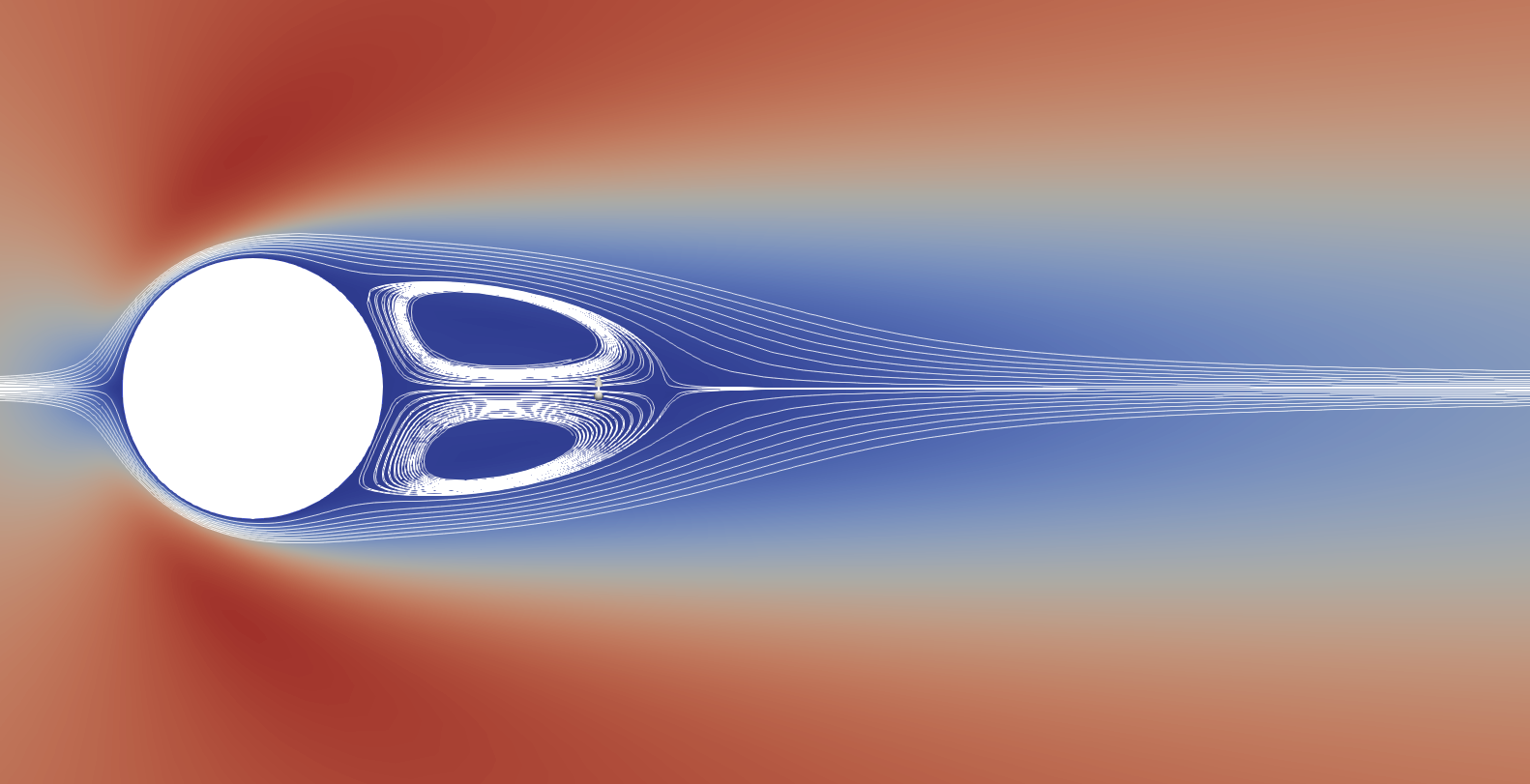}\\
\end{tabular}
\caption{Variation of the power index $m$. From top to bottom $m=0.4, 0.6$, the left column correspond to $C_U = 10$ the right column to $C_U = 20$.}
\label{bulles}
\end{figure}

Figure~\ref{bulles} presentes some streamlines, near the obstacle, for the velocity obtained with the SRP$_{im}$ for both $C_U$ number, with $m=0.4$ and $0.6$. This illustrates the variations of the wake lenght and overall behaviour of the method which is in agreement with the literature of the subject.

\begin{table}[!ht]
\centering
\begin{tabular}{cccccccccc}
\hline
& & & \multicolumn{3}{c}{$C_D$} & &\multicolumn{3}{c}{$L$}\\ 
\hline 
$C_U$ & $m$ &$\quad$& SS-Mix & SRP$_{im}$ & $\quad \Delta_{C_D}(\%)$ &$\quad$ & SS-Mix & SRP$_{im}$ & $\quad\Delta_{L}(\%)$\\ 
\hline   
\multirow{8}{*}{10}  

%& 0.3 && 1.1045  &  1.1060  &  0.14  &&  1.7145  &  1.7186  &  0.24\\
& 1.0  &&  2.7517  &  2.7537  &  0.07  &&  0.2364  &  0.2376  &  0.52 \\
& 0.9  &&  2.4697  &  2.4721  &  0.10  &&  0.3177  &  0.3205  &  0.87 \\
& 0.8  &&  2.1955  &  2.1963  &  0.04  &&  0.4280  &  0.4298  &  0.42 \\
& 0.7  &&  1.9292  &  1.9289  &  0.01  &&  0.5808  &  0.5816  &  0.12 \\
& 0.6  &&  1.6747  &  1.6735  &  0.07  &&  0.7955  &  0.7946  &  0.11 \\
& 0.5  &&  1.4348  &  1.4334  &  0.10  &&  1.1055  &  1.1031  &  0.22 \\
& 0.4  &&  1.2140  &  1.2150  &  0.08  &&  1.5617  &  1.5642  &  0.16 \\
\hline
\multirow{8}{*}{20}  
%& 0.3 &  &  0.8958  &  0.8971  &  0.15  &  &  3.36219  &  3.3629  &  0.02\\
& 1.0  &&  2.7517  &  2.7537  &  0.07  &&  0.2364  &  0.2376  &  0.52 \\
& 0.9  &&  2.3920  &  2.3908  &  0.05  &&  0.3703  &  0.3710  &  0.19 \\
& 0.8  &&  2.0524  &  2.0528  &  0.02  &&  0.5627  &  0.5639  &  0.22 \\
& 0.7  &&  1.7395  &  1.7404  &  0.05  &&  0.8437  &  0.8447  &  0.11 \\
& 0.6  &&  1.4583  &  1.4588  &  0.04  &&  1.2652  &  1.2662  &  0.08 \\
& 0.5  &&  1.2123  &  1.2118  &  0.04  &&  1.9086  &  1.9077  &  0.05 \\
& 0.4  &&  1.0042  &  1.0045  &  0.03  &&  2.8702  &  2.8709  &  0.02 \\
\hline
\end{tabular}
\caption{Drag coefficient ($C_D$) and wake lenght ($L$) at Reynolds $Re=10$, Carreau number $C_U=10, 20$ and for a power index ($m$) ranging from $0.4$ to $1.0$. Results for the steady states mixed (denoted SS-Mix) and unsteady SRP$_{im}$ schemes.}
\label{tabloCU10}
\end{table}
Table~\ref{tabloCU10} presents the drag coefficient and wake (or re-circulation) lenght for both values of the Carreau number ($C_U = 10$ and $20$) and a power index ranging from $0.4$ to $1.0$. Relative variations of the results from the steady state mixed method and the SRP$_{im}$ are also presented. In both cases these variations are roughly less than $1\%$, which can be attributed to the convergence criteria to reach the steady state, confirming the good behaviour of the SRP$_{im}$ method even for local quantities of interests such as the drag coefficient. Although not presented here, the results of Table~\ref{tabloCU10} for $C_U=10$ and $m=0.4, 0.6, 0.8, 1.0$ are in good agreement with the results presented in~\cite{LasPraGia2012} and~\cite{Pan2016} with less then $2.5\%$ of variation for the drag coefficient.

\section{Conclusions}
We proposed an original projection method for the numerical simulation of generalized Newtonian flow.  This method takes into account the explicit dependence of the viscosity upon velocity.  Based on this projection and the finite element method,  a  family  of  numerical  scheme using a second order time approximation was proposed.

Using a manufactured solution and the simulation of the flow past a cylinder, we illustrate the validity of four representatives of the family of numerical methods. As expected, explicit methods were less expensive from a computational point of view, however they exhibit the usual weakness of such methods (mainly weak order of precision). As for the implicit versions, the SRP$_{im}$ can be viewed as a generalized rotational projections. The behaviour of both IP$_{im}$ and SRP$_{im}$ was satisfactory. We underline that the SRP exhibit, as expected from~\cite{DetYak2017}, an improved accuracy in approximating the pressure when compared to the more intuitive incremental approach.

Finally, as in~\cite{DetJenYak2014}, these shear rate projections could 
be used to produce new \textit{coupled projection} methods applicable to generic  convection-diffusion phenomenon coupled with Newtonian/Non Newtonian flows (natural convection, coextrusion, etc.). Following the analysis in \cite{DetJenYak2014} these methods should offer important gain in performance and accuracy. 

\section*{References}
\biboptions{sort&compress}
\bibliographystyle{model1-num-names}
\bibliography{banquebiblio}

\begin{thebibliography}{36}
\expandafter\ifx\csname natexlab\endcsname\relax\def\natexlab#1{#1}\fi
\providecommand{\bibinfo}[2]{#2}
\ifx\xfnm\relax \def\xfnm[#1]{\unskip,\space#1}\fi
%Type = Article
\bibitem[{Guermond et~al.(2006)Guermond, Minev, and Shen}]{GueMinShe2006}
\bibinfo{author}{J.~L. Guermond}, \bibinfo{author}{P.~Minev},
  \bibinfo{author}{J.~Shen},
\newblock \bibinfo{title}{An overview of projection methods for incompressible
  flows},
\newblock \bibinfo{journal}{Comput. Methods Appl. Mech. Engrg.}
  \bibinfo{volume}{195} (\bibinfo{year}{2006}) \bibinfo{pages}{6011--6045}.
%Type = Article
\bibitem[{Timmermans et~al.(1996)Timmermans, Minev, and Van
  De~Vosse}]{TimMinVan1996}
\bibinfo{author}{L.~Timmermans}, \bibinfo{author}{P.~Minev},
  \bibinfo{author}{F.~Van De~Vosse},
\newblock \bibinfo{title}{An approximate projection scheme for incompressible
  flow using spectral elements},
\newblock \bibinfo{journal}{Int. J. Numer. Meth. Fluids.} \bibinfo{volume}{22}
  (\bibinfo{year}{1996}) \bibinfo{pages}{673 -- 688}.
%Type = Book
\bibitem[{Boffi et~al.(2013)Boffi, Brezzi, and Fortin}]{BofBreFor2013}
\bibinfo{author}{D.~Boffi}, \bibinfo{author}{F.~Brezzi},
  \bibinfo{author}{M.~Fortin}, \bibinfo{title}{Mixed Finite Element Methods and
  Applications}, volume~\bibinfo{volume}{44} of
  \textit{\bibinfo{series}{Springer Series in Computational Mathematics}},
  \bibinfo{publisher}{Springer}, \bibinfo{address}{Berlin, Heidelberg},
  \bibinfo{year}{2013}.
%Type = Book
\bibitem[{Girault and Raviart(1986)}]{GirRav1986}
\bibinfo{author}{V.~Girault}, \bibinfo{author}{P.-A. Raviart},
  \bibinfo{title}{Finite Element Methods for the {N}avier-{S}tokes Equations},
  volume~\bibinfo{volume}{5} of \textit{\bibinfo{series}{Springer Series in
  Computational Mathematics}}, \bibinfo{publisher}{Springer-Verlag},
  \bibinfo{address}{Berlin}, \bibinfo{year}{1986}. \bibinfo{note}{Theory and
  algorithms}.
%Type = Article
\bibitem[{Chorin(1968)}]{ChoAle1968}
\bibinfo{author}{A.~Chorin},
\newblock \bibinfo{title}{Numerical solution of the {N}avier--{S}tokes
  equations},
\newblock \bibinfo{journal}{Math. Comp.} \bibinfo{volume}{22}
  (\bibinfo{year}{1968}) \bibinfo{pages}{745--762}.
%Type = Article
\bibitem[{Chorin(1969)}]{ChoAle1969}
\bibinfo{author}{A.~Chorin},
\newblock \bibinfo{title}{On the convergence of discrete approximations to the
  {N}avier--{S}tokes equations},
\newblock \bibinfo{journal}{Math. Comp.} \bibinfo{volume}{23}
  (\bibinfo{year}{1969}) \bibinfo{pages}{341--353}.
%Type = Book
\bibitem[{Temam(1977)}]{Tem1977}
\bibinfo{author}{R.~Temam}, \bibinfo{title}{Navier--{S}tokes equations. Theory
  and Numerical Analysis}, volume~\bibinfo{volume}{2} of
  \textit{\bibinfo{series}{Studies in Mathematics and its Applications}},
  \bibinfo{publisher}{North-Holland Publishing}, \bibinfo{year}{1977}.
%Type = Incollection
\bibitem[{Rannacher(1992)}]{Ran1992}
\bibinfo{author}{R.~Rannacher},
\newblock \bibinfo{title}{On {C}horin's projection method for the
  incompressible {N}avier-{S}tokes equations},
\newblock in: \bibinfo{editor}{J.~Heywood}, \bibinfo{editor}{K.~Masuda},
  \bibinfo{editor}{R.~Rautmann}, \bibinfo{editor}{V.~Solonnikov} (Eds.),
  \bibinfo{booktitle}{The Navier-Stokes Equations II — Theory and Numerical
  Methods}, volume \bibinfo{volume}{1530} of \textit{\bibinfo{series}{Lecture
  Notes in Mathematics}}, \bibinfo{publisher}{Springer Berlin Heidelberg},
  \bibinfo{year}{1992}, pp. \bibinfo{pages}{167--183}.
%Type = Article
\bibitem[{Shen(1992)}]{She1992}
\bibinfo{author}{J.~Shen},
\newblock \bibinfo{title}{On error estimates of projection methods for
  {N}avier-{S}tokes equations: first-order schemes},
\newblock \bibinfo{journal}{SIAM J. Numer. Anal.} \bibinfo{volume}{29}
  (\bibinfo{year}{1992}) \bibinfo{pages}{57--77}.
%Type = Article
\bibitem[{Bell et~al.(1989)Bell, Colella, and Glaz}]{BelColGla1989}
\bibinfo{author}{J.~Bell}, \bibinfo{author}{P.~Colella},
  \bibinfo{author}{H.~Glaz},
\newblock \bibinfo{title}{A second-order projection method for the
  incompressible {N}avier-{S}tokes equations},
\newblock \bibinfo{journal}{J. of Comp. Physics} \bibinfo{volume}{85}
  (\bibinfo{year}{1989}) \bibinfo{pages}{257 -- 283}.
%Type = Article
\bibitem[{Shen(1996)}]{She1996}
\bibinfo{author}{J.~Shen},
\newblock \bibinfo{title}{On error estimates of the projection methods for the
  {N}avier-{S}tokes equations: second-order schemes},
\newblock \bibinfo{journal}{Math. Comp.} \bibinfo{volume}{65}
  (\bibinfo{year}{1996}) \bibinfo{pages}{1039--1065}.
%Type = Article
\bibitem[{Guermond(1997)}]{Gue1997}
\bibinfo{author}{J.-L. Guermond},
\newblock \bibinfo{title}{Un r\'esultat de convergence d'ordre deux pour
  l'approximation des \'equations de {N}avier-{S}tokes par projection
  incr\'ementale},
\newblock \bibinfo{journal}{C. R. Acad. Sci. Paris S\'er. I Math.}
  \bibinfo{volume}{325} (\bibinfo{year}{1997}) \bibinfo{pages}{1329--1332}.
%Type = Article
\bibitem[{Guermond and Quartapelle(1998)}]{GueQua1998}
\bibinfo{author}{J.-L. Guermond}, \bibinfo{author}{L.~Quartapelle},
\newblock \bibinfo{title}{On stability and convergence of projection methods
  based on pressure {P}oisson equation},
\newblock \bibinfo{journal}{Internat. J. Numer. Methods Fluids}
  \bibinfo{volume}{26} (\bibinfo{year}{1998}) \bibinfo{pages}{1039--1053}.
%Type = Article
\bibitem[{Goda(1979)}]{God1979}
\bibinfo{author}{K.~Goda},
\newblock \bibinfo{title}{A multistep technique with implicit difference
  schemes for calculating two- or three-dimensional cavity flows},
\newblock \bibinfo{journal}{J. Comput. Phys} \bibinfo{volume}{30}
  (\bibinfo{year}{1979}) \bibinfo{pages}{76--95}.
%Type = Article
\bibitem[{van Kan(1986)}]{Van1986}
\bibinfo{author}{J.~van Kan},
\newblock \bibinfo{title}{A second-order accurate pressure-correction scheme
  for viscous incompressible flow},
\newblock \bibinfo{journal}{SIAM J. Sci. Statist. Comput.} \bibinfo{volume}{7}
  (\bibinfo{year}{1986}) \bibinfo{pages}{870--891}.
%Type = Article
\bibitem[{Deteix and Yakoubi(2018)}]{DetYak2017}
\bibinfo{author}{J.~Deteix}, \bibinfo{author}{D.~Yakoubi},
\newblock \bibinfo{title}{Improving the pressure accuracy in a projection
  scheme for incompressible fluids with variable viscosity},
\newblock \bibinfo{journal}{Applied Mathematics Letters} \bibinfo{volume}{79}
  (\bibinfo{year}{2018}) \bibinfo{pages}{111 -- 117}.
%Type = Book
\bibitem[{Owens and Phillips(2002)}]{OwePhi2002}
\bibinfo{author}{R.~G. Owens}, \bibinfo{author}{T.~N. Phillips},
  \bibinfo{title}{Computational rheology}, volume~\bibinfo{volume}{14},
  \bibinfo{publisher}{World Scientific}, \bibinfo{year}{2002}.
%Type = Book
\bibitem[{Irgens(2013)}]{Irg2013}
\bibinfo{author}{F.~Irgens}, \bibinfo{title}{Rheology and {N}on-{N}ewtonian
  Fluids}, \bibinfo{publisher}{Springer International Publishing},
  \bibinfo{year}{2013}.
%Type = Article
\bibitem[{Bae and Wolf(2016)}]{BaeWol2016}
\bibinfo{author}{H.-O. Bae}, \bibinfo{author}{J.~Wolf},
\newblock \bibinfo{title}{Sufficient conditions for local regularity to the
  generalized newtonian fluid with shear thinning viscosity},
\newblock \bibinfo{journal}{Zeitschrift f{\"u}r angewandte Mathematik und
  Physik} \bibinfo{volume}{68} (\bibinfo{year}{2016}) \bibinfo{pages}{7}.
%Type = Article
\bibitem[{Tan and Zhou(2017)}]{TanZho2017}
\bibinfo{author}{Z.~Tan}, \bibinfo{author}{J.~Zhou},
\newblock \bibinfo{title}{Partial regularity of a certain class of
  non-newtonian fluids},
\newblock \bibinfo{journal}{Journal of Mathematical Analysis and Applications}
  \bibinfo{volume}{455} (\bibinfo{year}{2017}) \bibinfo{pages}{1529 -- 1558}.
%Type = Article
\bibitem[{Dreyfuss and Hungerb\"uhler(2004)}]{DreHun2004}
\bibinfo{author}{P.~Dreyfuss}, \bibinfo{author}{N.~Hungerb\"uhler},
\newblock \bibinfo{title}{Results on a {N}avier-{S}tokes system with
  applications to electrorheological fluid flow},
\newblock \bibinfo{journal}{Int. J. Pure Appl. Math.} \bibinfo{volume}{14}
  (\bibinfo{year}{2004}) \bibinfo{pages}{241--271}.
%Type = Article
\bibitem[{Diening et~al.(2010)Diening, R{\u{u}}{\v{z}}i{\v{c}}ka, and
  Wolf}]{DieRuzWol2010}
\bibinfo{author}{L.~Diening}, \bibinfo{author}{M.~R{\u{u}}{\v{z}}i{\v{c}}ka},
  \bibinfo{author}{J.~Wolf},
\newblock \bibinfo{title}{Existence of weak solutions for unsteady motions of
  generalized {N}ewtonian fluids},
\newblock \bibinfo{journal}{Ann. Sc. Norm. Super. Pisa Cl. Sci. (5)}
  \bibinfo{volume}{9} (\bibinfo{year}{2010}) \bibinfo{pages}{1--46}.
%Type = Article
\bibitem[{Bae(2015)}]{Bae2015}
\bibinfo{author}{H.-O. Bae},
\newblock \bibinfo{title}{Regularity criterion for generalized newtonian fluids
  in bounded domains},
\newblock \bibinfo{journal}{Journal of Mathematical Analysis and Applications}
  \bibinfo{volume}{421} (\bibinfo{year}{2015}) \bibinfo{pages}{489 -- 500}.
%Type = Article
\bibitem[{Hecht(2012)}]{freefem}
\bibinfo{author}{F.~Hecht},
\newblock \bibinfo{title}{New development in freefem++},
\newblock \bibinfo{journal}{J. Numer. Math.} \bibinfo{volume}{20}
  (\bibinfo{year}{2012}) \bibinfo{pages}{251--265}.
%Type = Book
\bibitem[{Adams and Fournier(2003)}]{Ada2003}
\bibinfo{author}{R.~A. Adams}, \bibinfo{author}{J.~J.~F. Fournier},
  \bibinfo{title}{Sobolev spaces}, volume \bibinfo{volume}{140},
  \bibinfo{publisher}{Pure and Applied Mathematics}, \bibinfo{address}{Academic
  Press, New York, London}, \bibinfo{year}{2003}.
%Type = Inbook
\bibitem[{Galdi et~al.(2000)Galdi, Heywood, and Rannacher}]{Ran2000}
\bibinfo{editor}{G.~P. Galdi}, \bibinfo{editor}{J.~G. Heywood},
  \bibinfo{editor}{R.~Rannacher} (Eds.), \bibinfo{title}{Finite Element Methods
  for the Incompressible {N}avier--{S}tokes Equations},
  \bibinfo{publisher}{Birkh{\"a}user Basel}, \bibinfo{address}{Basel}, pp.
  \bibinfo{pages}{191--293}.
%Type = Article
\bibitem[{Guermond and Shen(2004)}]{GueShe2004}
\bibinfo{author}{J.~L. Guermond}, \bibinfo{author}{J.~Shen},
\newblock \bibinfo{title}{On the error estimates for the rotational
  pressure-correction projection methods.},
\newblock \bibinfo{journal}{Math. Comput.} \bibinfo{volume}{73}
  (\bibinfo{year}{2004}) \bibinfo{pages}{1719--1737}.
%Type = Book
\bibitem[{Boyer and Fabrie(2006)}]{BoyFab2006}
\bibinfo{author}{F.~Boyer}, \bibinfo{author}{P.~Fabrie},
  \bibinfo{title}{\'El\'ements d'analyse pour l'\'etude de quelques mod\`eles
  d'\'ecoulements de fluides visqueux incompressibles},
  volume~\bibinfo{volume}{52} of \textit{\bibinfo{series}{Math\'ematiques \&
  Applications (Berlin) [Mathematics \& Applications]}},
  \bibinfo{publisher}{Springer-Verlag, Berlin}, \bibinfo{year}{2006}.
%Type = Article
\bibitem[{Guermond et~al.(2005)Guermond, Minev, and Shen}]{GueMinShe2005}
\bibinfo{author}{J.~L. Guermond}, \bibinfo{author}{P.~Minev},
  \bibinfo{author}{J.~Shen},
\newblock \bibinfo{title}{Error analysis of pressure-correction schemes for the
  time-dependent stokes equations with open boundary conditions},
\newblock \bibinfo{journal}{SIAM Journal on Numerical Analysis}
  \bibinfo{volume}{43} (\bibinfo{year}{2005}) \bibinfo{pages}{239--258}.
%Type = Book
\bibitem[{Bathe(1996)}]{Bat1996}
\bibinfo{author}{K.~Bathe}, \bibinfo{title}{Finite Element Procedures},
  \bibinfo{publisher}{Prentice-Hall}, \bibinfo{address}{New Jersey},
  \bibinfo{year}{1996}.
%Type = Book
\bibitem[{Ciarlet and Luneville(2009)}]{CiaLun2009}
\bibinfo{author}{P.~Ciarlet}, \bibinfo{author}{E.~Luneville},
  \bibinfo{title}{La m{\'e}thode des {\'e}l{\'e}ments finis: de la th{\'e}orie
  {\`a} la pratique. Concepts g{\'e}n{\'e}raux. I}, Cours (ENSTA),
  \bibinfo{publisher}{Les Presses de l'ENSTA}, \bibinfo{year}{2009}.
%Type = Article
\bibitem[{Lashgari et~al.(2012)Lashgari, Pralits, Giannetti, and
  Brandt}]{LasPraGia2012}
\bibinfo{author}{I.~Lashgari}, \bibinfo{author}{J.~O. Pralits},
  \bibinfo{author}{F.~Giannetti}, \bibinfo{author}{L.~Brandt},
\newblock \bibinfo{title}{First instability of the flow of shear-thinning and
  shear-thickening fluids past a circular cylinder},
\newblock \bibinfo{journal}{Journal of Fluid Mechanics} \bibinfo{volume}{701}
  (\bibinfo{year}{2012}) \bibinfo{pages}{201–227}.
%Type = Article
\bibitem[{Patnana et~al.(2009)Patnana, Bharti, and Chhabra}]{PatBhaChh2009}
\bibinfo{author}{V.~Patnana}, \bibinfo{author}{R.~Bharti},
  \bibinfo{author}{R.~Chhabra},
\newblock \bibinfo{title}{Two-dimensional unsteady flow of power-law fluids
  over a cylinder},
\newblock \bibinfo{journal}{Chemical Engineering Science} \bibinfo{volume}{64}
  (\bibinfo{year}{2009}) \bibinfo{pages}{2978 -- 2999}.
%Type = Article
\bibitem[{Pantokratoras(2016)}]{Pan2016}
\bibinfo{author}{A.~Pantokratoras},
\newblock \bibinfo{title}{Steady flow of a non-newtonian carreau fluid across
  an unconfined circular cylinder},
\newblock \bibinfo{journal}{Meccanica} \bibinfo{volume}{51}
  (\bibinfo{year}{2016}) \bibinfo{pages}{1007--1016}.
%Type = Article
\bibitem[{Mossaz et~al.(2010)Mossaz, Jay, and Magnin}]{MosJayMag2010}
\bibinfo{author}{S.~Mossaz}, \bibinfo{author}{P.~Jay},
  \bibinfo{author}{A.~Magnin},
\newblock \bibinfo{title}{Criteria for the appearance of recirculating and
  non-stationary regimes behind a cylinder in a viscoplastic fluid},
\newblock \bibinfo{journal}{Journal of Non-Newtonian Fluid Mechanics}
  \bibinfo{volume}{165} (\bibinfo{year}{2010}) \bibinfo{pages}{1525 -- 1535}.
%Type = Article
\bibitem[{Deteix et~al.(2014)Deteix, Jendoubi, and Yakoubi}]{DetJenYak2014}
\bibinfo{author}{J.~Deteix}, \bibinfo{author}{A.~Jendoubi},
  \bibinfo{author}{D.~Yakoubi},
\newblock \bibinfo{title}{A coupled prediction scheme for solving the
  {N}avier--{S}tokes and heat equations},
\newblock \bibinfo{journal}{SIAM Journal of Numerical Analysis}
  \bibinfo{volume}{52} (\bibinfo{year}{2014}) \bibinfo{pages}{2415--2439}.

\end{thebibliography}
% \bibliography{biblio/banquebiblio}
 %     \bibliography{banquebiblio.bib}
 \end{document}